\documentclass[10pt]{amsart}

\usepackage{amsmath}
\usepackage{amssymb}
\usepackage{enumerate}
\usepackage{amsbsy}
\usepackage{amsfonts}
\usepackage{color}

\usepackage[T1]{fontenc}
\usepackage{lmodern,microtype}
\usepackage[colorlinks=true,linkcolor=blue,citecolor=blue,urlcolor=blue]{hyperref}

\newtheorem{theorem}{Theorem}[section]
\newtheorem{lemma}[theorem]{Lemma}
\newtheorem{proposition}[theorem]{Proposition}
\newtheorem{corollary}[theorem]{Corollary}

\theoremstyle{definition}\newtheorem{assumption}[theorem]{Assumption}
\theoremstyle{remark}\newtheorem{remark}[theorem]{Remark}
\numberwithin{equation}{section}

\newcommand{\R}{\mathbb R}
\newcommand{\D}{\mathcal D}

\newcommand{\Op}{\operatorname{Op}}

\allowdisplaybreaks

\begin{document}
	
	\title[Dirac scattering for decaying spacetime metrics]{Scattering for cubic and Yukawa--Hartree Dirac equations wiith polynomially decaying spacetime metrics}

	\author[Y. Cho]{Yonggeun Cho}
	\address{ Department of Mathematics, and Institute of Pure and Applied Mathematics, Jeonbuk National University, Jeonju 54896, Republic of Korea}
	\email{changocho@jbnu.ac.kr}

\begin{abstract}
We study massive Dirac equations in three space dimensions with small,
time-dependent perturbations of the principal coefficients and decaying
zeroth-order matrix terms. The principal perturbation and its spatial
derivatives decay polynomially at rate $\rho>7/2$, and its first time
derivative has one additional inverse spatial power. We prove small-data
global well-posedness and scattering for local cubic nonlinearities in
$H^s$, $s>1+5/(2(\rho-1))$, and for Yukawa--Hartree nonlinearities in
$H^s$, $s>1/2+5/(2(\rho-1))$. The argument retains the coupling between
the metric-adapted positive and negative branches. Its main ingredients
are a time-independent modified energy, a weighted estimate controlling
the full massive low-frequency component, and frequency-dependent
spatial cutoffs. The operator class contains the geometric Dirac
Hamiltonian of a synchronous spacetime metric with its full spin
connection. No second time derivative, temporal integrability, or
asymptotic stationarity is assumed. Scattering is formulated relative
to the exact nonautonomous linear evolution.
\end{abstract}


\maketitle

\section{Introduction}
\label{sec:introduction}

The long-time analysis of a nonlinear Dirac equation with variable
principal coefficients involves both dispersive and geometric issues.
The perturbation of the flat Hamiltonian is of first order and cannot
be treated as a bounded forcing term in a direct global Duhamel
argument. In addition, the eigenspaces of the principal symbol depend
on spacetime. Applying their quantizations creates commutator and
composition terms, and the positive and negative branches remain
coupled. For a nonstationary metric, these terms contain the temporal
variation of the frame and the full spin connection.

In this paper we establish a small-data scattering theory under
explicit polynomial spatial bounds. The principal point is a global
weighted energy estimate for the full coupled Hamiltonian. It is based
on a fixed positive commutator and therefore does not accumulate the
coefficient variation through a global Gronwall inequality. The
weighted estimate absorbs the zeroth-order coupling directly, including
the low-frequency terms associated with the mass. We then combine it
with free massive estimates outside a frequency-dependent spatial
region. A local wave parametrix controls the remaining interior
contribution for the local cubic equation. The resulting regularity
loss is expressed explicitly in terms of the spatial decay rate.

\subsection{The operator class}
Throughout the paper $x\in\mathbb R^3$, $t\in\mathbb R$, and the mass
$M$ is strictly positive. Put $D_t=-i\partial_t$ and $D_j=-i\partial_j$.
The constant Hermitian $4\times4$ Dirac matrices satisfy
\[
 \alpha^a\alpha^b+\alpha^b\alpha^a=2\delta_{ab}I_4,
 \qquad \alpha^a\beta+\beta\alpha^a=0,
 \qquad \beta^2=I_4.
\]
Let $G(t,x)=I_3+h(t,x)$ be real symmetric and positive definite,
and choose its positive symmetric square root $e=G^{1/2}$, with
$g^{ij}=\sum_a e_a^{\ i}e_a^{\ j}$. We consider
\begin{equation}\label{main:hamiltonian}
 B(t)=\frac12\sum_{a,j=1}^3\alpha^a
       \bigl(e_a^{\ j}(t,x)D_j+D_je_a^{\ j}(t,x)\bigr)
       +M\beta+C(t,x),
\end{equation}
where $C$ is a complex $4\times4$ matrix. In the second summand of the
symmetrized expression, $D_je_a^{\ j}$ denotes operator composition.
The matrix $C$ need not be Hermitian. This formulation includes the
zeroth-order terms resulting from the spin connection and the change
to Weyl quantization.

For $\rho>0$ and an integer $J\ge1$, define
\begin{align*}
 \mathfrak g_{\rho,J}(h)
 &=\sum_{|\gamma|\le J}\sup_{t,x}
 \left(\langle x\rangle^{\rho+|\gamma|}|\partial_x^\gamma h|
 +\langle x\rangle^{\rho+1+|\gamma|}
                         |\partial_t\partial_x^\gamma h|\right),\\
 \mathfrak c_{\rho+1,J}(C)
 &=\sum_{|\gamma|\le J}\sup_{t,x}
          \langle x\rangle^{\rho+1+|\gamma|}
                         |\partial_x^\gamma C|.
\end{align*}
Matrix norms are any fixed equivalent finite-dimensional norms.

\begin{assumption}\label{ass:coefficients}
The coefficient $h$ is $C^1$ in time with values in $C_b^J$, and
$C$ is continuous in time with values in $C_b^J$. Their bounds satisfy
\begin{equation}\label{main:coeff}
 \mathfrak g_{\rho,J}(h)+\mathfrak c_{\rho+1,J}(C)\le\varepsilon.
\end{equation}
The integer $J$ is chosen sufficiently large for the fixed Sobolev
indices under consideration and $\varepsilon$ is small.
\end{assumption}

In particular, $G$ is uniformly elliptic if $\varepsilon$ is small.
There is no assumption on $\partial_tC$ or on a second time derivative
of $h$. All estimates depend on finitely many of the displayed
coefficient bounds. The proof first applies to smooth coefficients;
the finite regularity assertion follows by approximation with a
fixed derivative margin.

The nonlinear Cauchy problem is
\begin{equation}\label{main:equation}
 (D_t+B(t))\psi=\lambda\mathcal N(\psi),\qquad
 \psi(0)=\psi_0\in H^s(\mathbb R^3;\mathbb C^4),
 \qquad \lambda\in\mathbb R.
\end{equation}
Thus $i\partial_t\psi=B(t)\psi-\lambda\mathcal N(\psi)$.
We allow either sign of $\lambda$ and fixed finite sums of the forms
\begin{align}
 \mathcal N_C(\psi)&=(\psi^\dagger C_1\psi)C_2\psi,
             \label{main:cubic}\\
 \mathcal N_H(\psi)&=\bigl(V_b*(\psi^\dagger C_1\psi)\bigr)C_2\psi,
 \qquad V_b(x)=\frac{e^{-b|x|}}{|x|},\quad b>0,
             \label{main:hartree}
\end{align}
where $C_1,C_2$ are fixed complex matrices. In particular, scalar-density
interactions and scalar bilinear interactions involving $\beta$ are
included. The convolution in \eqref{main:hartree} is taken with respect
to Lebesgue measure in the specified spatial coordinates.

\subsection{Main results}
Let $U_D(t,r)$ denote the exact evolution of $(D_t+B(t))u=0$, with
$U_D(r,r)=I$. It obeys the cocycle identity and
$U_D(r,t)=U_D(t,r)^{-1}$. No identification of the inverse with the
Lebesgue $L^2$ adjoint is needed.

\begin{theorem}[Local cubic interactions]\label{main:theorem}
Fix $M>0$, $\rho>7/2$, and
\begin{equation}\label{main:cubic-range}
 s>1+\frac{5}{2(\rho-1)}.
\end{equation}
There exist an integer $J_0$, a coefficient size $\varepsilon_*>0$,
and a data size $\delta_*>0$ such that the following holds. Suppose
Assumption~\ref{ass:coefficients} holds with $J\ge J_0$ and
$\varepsilon\le\varepsilon_*$. If $\|\psi_0\|_{H^s}\le\delta_*$,
then \eqref{main:equation} with \eqref{main:cubic} has a global solution
\[
 \psi\in C(\mathbb R;H^s)\cap L_t^\infty H_x^s
                         \cap L_t^2L_x^\infty,
 \qquad
 \|\psi\|_{L_t^\infty H_x^s}+\|\psi\|_{L_t^2L_x^\infty}
 \le C_s\|\psi_0\|_{H^s}.
\]
It is unique in this class and depends locally Lipschitz continuously
on the data. There are $\psi_\pm\in H^s$ such that
\begin{equation}\label{main:scatter}
 \lim_{t\to\pm\infty}
       \|\psi(t)-U_D(t,0)\psi_\pm\|_{H^s}=0.
\end{equation}
Conversely, every sufficiently small prescribed profile at either
time infinity is attained by a unique small solution in this class.
\end{theorem}

For the Hartree equation we use a quadratic spacetime norm $Q_s$,
defined by frequency cubes in Section~\ref{dirac:fourier}, and set
$W^s=L_t^\infty H_x^s\cap\{Q_s<\infty\}$ with the sum norm.
Its exponent is chosen as in \eqref{dirac:hartree-parameters}.

\begin{theorem}[Yukawa--Hartree interactions]\label{main:hartree-theorem}
Fix $M,b>0$, $\rho>7/2$, and
\begin{equation}\label{main:hartree-range}
 s>\frac12+\frac{5}{2(\rho-1)}.
\end{equation}
For sufficiently large $J$ and sufficiently small coefficient and
data sizes, as in Theorem~\ref{main:theorem}, problem
\eqref{main:equation} with \eqref{main:hartree} has a global solution
\[
 \psi\in C(\mathbb R;H^s)\cap W^s,
 \qquad \|\psi\|_{W^s}\le C_s\|\psi_0\|_{H^s}.
\]
Uniqueness in this class, local Lipschitz dependence, the scattering
limits \eqref{main:scatter}, and small wave operators hold.
\end{theorem}

The constants in these statements may depend on the fixed mass,
decay exponent, Sobolev index, coupling, and nonlinearity matrices;
in the Hartree case they may also depend on $b$. The endpoints are
not included. The conclusions do not require a null structure of
the nonlinearity.

The common linear estimate is useful independently of the nonlinear
applications. Its coefficient decay requirement is weaker than that
needed for the dispersive estimates.

\begin{theorem}[Uniform weighted energy]\label{main:linear-theorem}
Let $\rho>2$, $s\ge0$, and
\[
 0<\nu<\min\{1,\rho-2\},\qquad
 a=\frac{3+\nu}{2},\qquad w=\langle x\rangle^{-a}.
\]
Under Assumption~\ref{ass:coefficients}, with sufficiently large $J$
and sufficiently small $\varepsilon$, the full evolution satisfies
\begin{equation}\label{main:linear}
 \sup_{t,r}\|U_D(t,r)\|_{H^s\to H^s}\le C_s,
 \qquad
 \|wU_D(\cdot,r)f\|_{L_t^2H_x^s}\le C_s\|f\|_{H^s}.
\end{equation}
More generally, on either time orientation of an interval $I$ with
initial endpoint $t_0$, solutions of $(D_t+B)u=F_1+F_2$ satisfy
\begin{equation}\label{main:forced-linear}
 \begin{split}
 \|u\|_{L^\infty(I;H^s)}+\|wu\|_{L^2(I;H^s)}
 \le C_s\bigl(&\|u(t_0)\|_{H^s}+\|F_1\|_{L^1(I;H^s)}\\
              &+\|w^{-1}F_2\|_{L^2(I;H^s)}\bigr).
 \end{split}
\end{equation}
The constant is independent of the interval and its initial time.
\end{theorem}

\begin{corollary}[Synchronous geometric backgrounds]\label{main:geometric}
Let $\mathbf g=-dt^2+\gamma_{ij}(t,x)dx^idx^j$, with
$G=\gamma^{-1}=I+h$, and suppose
$\mathfrak g_{\rho,J+1}(h)$ is sufficiently small. In the instantaneous
symmetric orthonormal frame, the full geometric Dirac Hamiltonian
obeys the conclusions above in the corresponding ranges of $\rho,s$.
The nonlinearities are \eqref{main:cubic} and \eqref{main:hartree}
in that frame and those coordinates.
\end{corollary}
The verification that the full spin connection fits
\eqref{main:hamiltonian} is given in Section~\ref{sec:geometry}.
The extra spatial derivative in the corollary accounts for the
spatial connection coefficients. Neither the corollary nor the
analytic theorems require a bound on $\partial_t^2h$.

\subsection{Examples and the decay--regularity tradeoff}
The metric may vary persistently in time. For a fixed real symmetric
matrix $H_*$, put
\[
 G(t,x)=I_3+\eta a(t)\langle x\rangle^{-\rho-1}H_*.
\]
If $a,a'$ are bounded and $\eta$ is sufficiently small, the coefficient
hypotheses hold with any prescribed finite number of spatial
derivatives. The choice $a(t)=\sin t$ is periodic and has no limiting
metric. A further example is
\[
 a(t)=\sin t+\int_0^t\sin(r^2)\,dr.
\]
The oscillatory integral is bounded, $a'$ is bounded, and
$a''(t)=-\sin t+2t\cos(t^2)$ is unbounded. Thus the hypotheses
also allow metrics without a uniform second-time-derivative bound.
No compact spatial support is present in either example.

For $\rho=6$, the sufficient regularities are $s>3/2$ for the cubic
equation and $s>1$ for Hartree. For $\rho=11$, they are $s>5/4$
and $s>3/4$, respectively. Given any fixed $s_0>1$ in the cubic
case, choosing
\[
 \rho>\max\left\{\frac72,1+\frac{5}{2(s_0-1)}\right\}
\]
makes Theorem~\ref{main:theorem} available at each fixed $s\ge s_0$.
Given $s_0>1/2$ in the Hartree case, the analogous sufficient choice is
\[
 \rho>\max\left\{\frac72,1+\frac{5}{2s_0-1}\right\}.
\]
These thresholds describe the loss of the present argument and are
not asserted to be optimal. The number of coefficient derivatives
and the smallness constants may depend on the chosen $s$.

\subsection*{Related work and position of the results}

The small-data theory of nonlinear Dirac equations on Minkowski space is
substantially more precise than the Sobolev thresholds proved here.
Machihara, Nakanishi, and Ozawa~\cite{MNO2003} established a small-data
theory for the massive nonlinear Dirac equation in $H^s(\mathbb R^3)$,
$s>1$, using spacetime estimates for the Klein--Gordon evolution.
For the standard cubic Dirac nonlinearity, Bejenaru and Herr~\cite{BH2015}
proved global well-posedness and scattering at the critical regularity
$H^1(\mathbb R^3)$.  Their argument uses the particular algebraic structure
of the nonlinearity and refined estimates for the free evolution.
The present work addresses the effect of nonstationary principal
coefficients, with a loss that is quantified by the rate of spatial decay;
it does not improve these flat-space regularity results.

The flat Hartree theory also reaches lower regularities than those considered
here.  Herr and Tesfahun~\cite{HT} proved small-data scattering for
semi-relativistic Hartree equations, including the nonradial massive case
at $s>1/2$, and discussed a Dirac application.  For Dirac equations with
the standard null-structured Hartree interaction, Yang~\cite{Yang2019}
and Tesfahun~\cite{Tesfahun2020} obtained small-data scattering throughout
the subcritical range $s>0$.  Cho, Hong, and Lee~\cite{CHL2023} treated
Yukawa interactions at the critical $L^2$ regularity with additional
angular regularity.  These results use features of the flat evolution
and, in the latter works, the null structure of the interaction.  Our
Hartree estimate instead uses the smoothing of the Yukawa convolution
and the weighted spacetime bounds for the full variable-coefficient
Dirac evolution.  Throughout this paper the Yukawa interaction is the
specified convolution in the chosen Euclidean spatial coordinates;
no identification with a Green operator of the evolving spatial metric
is intended.

For linear Dirac equations with electromagnetic potentials,
D'Ancona and Fanelli~\cite{DF2008} proved global smoothing and Strichartz
estimates under suitable decay and smallness assumptions.  Variable
principal coefficients require a different treatment.  On static
asymptotically flat backgrounds, Cacciafesta and de Suzzoni~\cite{CdS2019}
proved local smoothing estimates, and Cacciafesta, de Suzzoni, and
Meng~\cite{CdSM2024} derived Strichartz estimates.  The metric in the
latter work has spatial part $h_{ij}(x)$ independent of time.  Our
coefficients may depend on both $t$ and $x$, and the estimates are
formulated uniformly on the whole time axis.

The relevant parametrix methods have an established history.
Smith~\cite{Smith1998} constructed local wave parametrices under
low-regularity coefficient assumptions, while Metcalfe and
Tataru~\cite{MT2012} developed global outgoing parametrices and
dispersive estimates for variable-coefficient wave equations, including
time-dependent metrics satisfying weak asymptotic flatness conditions.
Herr and Hong~\cite{HH2026} state a small-data scattering theorem for
the cubic massive Dirac equation throughout $s>1$ (Theorem~1.2).
In the third arXiv version, dated 14 July 2026, the reduction following
equations (2.12)--(2.13) contains a gap: the displayed summed equation
does not imply the two separate half-wave equations asserted there.
Remark~\ref{remark:HH-gap} gives the algebraic obstruction. We therefore
do not use that Dirac theorem or its decoupling as an established input.

Our purpose is to give an independent argument for the coupled Dirac
system under explicit polynomial spatial bounds.  We retain the
zeroth-order matrix terms, use the scalar principal symbol only to
organize the linear estimates, and establish scattering relative to the
full nonautonomous Dirac propagator.  The hypotheses require only first
time derivatives of the principal coefficients; the additional matrix
potential is continuous in time and is controlled through its spatial
derivatives.
These thresholds express a spatial decay--regularity tradeoff in this
proof and are not asserted to be optimal.  Compared with~\cite{HH2026},
the spatial assumptions are more restrictive and the Sobolev range is
narrower, whereas the time-regularity assumptions differ.  Accordingly,
the two sets of hypotheses are not ordered simply by inclusion.  No convergence of the coefficients as $t\to\pm\infty$
is assumed, and scattering to a prescribed Minkowski evolution is not
part of the conclusion.

\subsection{Outline of the proof}
The scalar branches are defined by the eigenvalues
$\pm q_g$, where $q_g=(M^2+g^{ij}\xi_i\xi_j)^{1/2}$. Their
quantized eigenprojections have exact sum $I_4$, but are not exact
operator projections. We retain the associated order-zero errors
and obtain an exact coupled system with remainders of spatial
order $-\rho-1$.

For the scalar Weyl generator $A=\Op^w(q_g)$, the perturbation
$A-(M^2-\Delta)^{1/2}$ has an explicit symbol of orders $(1,-\rho)$.
A fixed radial multiplier produces a commutator that controls the
full weighted Sobolev norm, including the massive low frequencies.
The associated energy is independent of time. Conjugating the
coupling by $w^{-1}$ on both sides gives an order-zero operator
with spatial order $2a-\rho-1<0$, so the matrix coupling is absorbed
in the weighted forcing estimate. This proves
Theorem~\ref{main:linear-theorem}.

The remaining step uses a cutoff at radius $R_N$ for each dyadic
frequency $N$. The boundary and metric-tail source terms cost
$R_N^{a+1}$ and $NR_N^{a+2-\rho}$, respectively. Their balance
explains the loss $(a+1)/(\rho-1)$, which approaches
$5/(2(\rho-1))$ as $a\downarrow3/2$. For the cubic equation, a
uniform local Strichartz estimate replaces the interior Bernstein
loss. We derive it by squaring only the scalar principal branches,
which introduces $\partial_tA$ but never differentiates the matrix
connection. For Hartree, frequency cubes and the two derivatives of
Yukawa smoothing yield the required trilinear estimate directly.
All nonlinear estimates are in physical spacetime norms, and the
final Duhamel and scattering arguments use the exact full evolution.

Section~\ref{dirac:reduction} develops the coupled reduction and
the geometric application. Sections~\ref{dirac:scalar-input} and
\ref{dirac:energy} prove the weighted energy estimates.
Section~\ref{dirac:fourier} contains the free and retarded estimates.
Sections~\ref{sec:gluing} and \ref{dirac:local-cubic} establish the
global resolution bounds. Section~\ref{sec:nonlinear} proves the
nonlinear theorems. Appendix~\ref{sec:atomic} records compatibility
with atomic spaces adapted to the nonautonomous evolution.

\section{Symbol calculus and the coupled Dirac reduction}
\label{dirac:reduction}

\subsection{Notation and weighted symbols}
All Sobolev spaces are inhomogeneous. We write $H^s$ for
$H^s(\mathbb R^3;\mathbb C^4)$, or its scalar counterpart as appropriate,
and use the Lebesgue spatial measure in all norms. The Fourier transform
is $\widehat f(\xi)=\int e^{-ix\cdot\xi}f(x)\,dx$.
The notation $A\lesssim B$ means $A\le CB$, with constants depending
only on the fixed parameters and finitely many specified coefficient
bounds. Constants are uniform in the position and length of a time
interval unless stated otherwise.

For a scalar or matrix symbol $a$, its Weyl quantization is
\[
 \Op^w(a)f(x)=(2\pi)^{-3}\iint e^{i(x-y)\cdot\xi}
               a\left(t,\frac{x+y}{2},\xi\right)f(y)\,dy\,d\xi.
\]
The class $SG^{m,\ell}$ consists of symbols satisfying
\begin{equation}\label{symbol:class}
 |\partial_x^\gamma\partial_\xi^\delta a(t,x,\xi)|
 \le C_{\gamma\delta}\langle x\rangle^{\ell-|\gamma|}
                         \langle\xi\rangle^{m-|\delta|}.
\end{equation}
The superscripts give the frequency and spatial orders, respectively.
For coefficients of finite smoothness, the notation means the finite
collection of seminorms required by the stated argument. The notation
$\varepsilon SG^{m,\ell}$ records an additional factor $\varepsilon$
in those seminorms. We use $\Psi^m$ for operators with uniformly
bounded symbols in $S^m_{1,0}$.

The composition and boundedness facts needed below follow from the
usual pseudodifferential calculus; see, for example,~\cite{Shubin2001}.
We record the orders to keep the role of spatial decay explicit.
If $a\in SG^{m_1,\ell_1}$ and $b\in SG^{m_2,\ell_2}$, then
\begin{equation}\label{symbol:composition}
 \Op^w(a)\Op^w(b)=\Op^w(ab)+\Op^w(r),
 \qquad r\in SG^{m_1+m_2-1,\ell_1+\ell_2-1}.
\end{equation}
The product $ab$ is a matrix product when appropriate. If it commutes
pointwise, the commutator has the order of $r$. An operator in
$SG^{0,\ell}$, $\ell\le0$, is bounded on $L^2$; conjugation by
massive Bessel multipliers gives its boundedness on any fixed $H^s$.
Multiplication by $\langle x\rangle^b$ changes the spatial order by
$b$ in these compositions.

For completeness, each term after $ab$ in the finite composition
expansion contains a frequency derivative paired with a spatial
derivative, giving the two order reductions in
\eqref{symbol:composition}. The oscillatory-integral remainder is
estimated after sufficiently many integrations by parts.
The inequalities
$\langle x\rangle^b\langle y\rangle^{-b}
\lesssim_b\langle x-y\rangle^{|b|}$ allow growing polynomial
weights in the same argument. Increasing the finite expansion order
makes the off-diagonal remainder integrable with these weights and
with the derivatives needed for the $L^2$ estimate. Thus each use
of the calculus below consumes only finitely many spatial derivatives;
no time derivative is involved in these composition rules.

\subsection{Metric-adapted half-wave components}
The Weyl symbol of the principal part of \eqref{main:hamiltonian} is
\[
 b_M(t,x,\xi)=\alpha^ae_a^{\ j}(t,x)\xi_j+M\beta,
 \qquad q_g(t,x,\xi)=\bigl(M^2+g^{ij}(t,x)\xi_i\xi_j\bigr)^{1/2}.
\]
The Clifford relations give the exact matrix identity
$b_M^2=q_g^2I_4$. Define
\begin{equation}\label{reduction:projectors}
 p_\sigma=\frac12\left(I_4+\sigma\frac{b_M}{q_g}\right),
 \qquad P_\sigma=\Op^w(p_\sigma),\qquad
 A=\Op^w(q_g),\qquad \sigma\in\{+1,-1\}.
\end{equation}
Pointwise, $p_\sigma^2=p_\sigma$,
$p_+p_-=0$, and $p_\sigma b_M=\sigma q_gp_\sigma$.
At the operator level we have the exact identity
\begin{equation}\label{reduction:sum}
 P_++P_-=I_4,
\end{equation}
whereas the idempotence and orthogonality identities acquire lower
order errors.

\begin{lemma}[Exact coupled equations]\label{reduction:lemma}
Under Assumption~\ref{ass:coefficients},
\[
 P_\sigma^2-P_\sigma,\ P_+P_-
       \in\varepsilon SG^{-1,-\rho-1}.
\]
If $(D_t+B)\psi=F$ and $u_\sigma=P_\sigma\psi$, then
\begin{equation}\label{dirac:coupled}
 (D_t+\sigma A)u_\sigma
 =P_\sigma F+R_\sigma(u_++u_-),
\end{equation}
where
\begin{equation}\label{dirac:remainder}
 R_\sigma=[D_t,P_\sigma]+\sigma AP_\sigma-P_\sigma B
       \in\varepsilon SG^{0,-\rho-1}.
\end{equation}
The identity \eqref{dirac:coupled} holds without discarding any
connection or projection terms.
\end{lemma}
\begin{proof}
Let $p_\sigma^0(\xi)$ be the flat symbols obtained by setting $e=I_3$.
The coefficient bounds and the positive mass give
\[
 p_\sigma-p_\sigma^0\in\varepsilon SG^{0,-\rho},\qquad
 \partial_tp_\sigma\in\varepsilon SG^{0,-\rho-1},\qquad
 \partial_xp_\sigma\in\varepsilon SG^{0,-\rho-1}.
\]
The projection defects follow from \eqref{symbol:composition}
and the pointwise identities. In the expression
$\sigma AP_\sigma-P_\sigma\Op^w(b_M)$ the pointwise order-one
terms cancel, leaving compositions with a spatial derivative of a
coefficient. The product $P_\sigma C$ has orders $(0,-\rho-1)$,
and $[D_t,P_\sigma]=-i\Op^w(\partial_tp_\sigma)$ has the same
orders. This proves \eqref{dirac:remainder}. Applying the operators
to the equation and using \eqref{reduction:sum} proves
\eqref{dirac:coupled}.
\end{proof}

The reduction is reversible with compatible data. Suppose
\eqref{dirac:coupled} is solved with $F=F(u_++u_-)$ and
$u_\sigma(r)=P_\sigma(r)f$. Put
$\psi=u_++u_-$ and $z=u_+-P_+\psi$. Since
\[
 R_++R_-=A(P_+-P_-)-B,
\]
direct summation and subtraction give
\begin{equation}\label{reduction:compatibility}
 (D_t+B)\psi=F(\psi)-2Az,
 \qquad D_tz=(2P_+-I_4)Az,
 \qquad z(r)=0.
\end{equation}
The last generator has Hermitian order-one principal symbol $b_M$.
Finite-interval energy uniqueness implies $z=0$. For mild solutions
the calculation is justified by spatial regularization and passage
to the limit. Thus the coupled system reconstructs the original
Dirac equation; exact idempotence of $P_\sigma$ is unnecessary.

\begin{remark}[A version-specific obstruction to the proposed decoupling]
\label{remark:HH-gap}
In this remark only, $\mathsf P_\pm$ denotes the flat-normalized
operators, in contrast to the metric-adapted $P_\pm$ above.
We describe the issue in~\cite[Section~2,
pp.~8--12]{HH2026} using notation that separates the flat scalar
multiplier from the full Dirac operator.  Write the homogeneous equation
in the form
\[
 (D_t+B(t))\psi=0,\qquad
 Q=(M^2+|D_x|^2)^{1/2},\qquad D_t=-i\partial_t,
\]
and define exactly as in the formula following (2.8) of that paper
\[
 T=Q^{-1}B,\qquad \mathsf P_\pm=\tfrac12(I\pm T),\qquad
 E=\tfrac14(T^2-I).
\]
These definitions imply the exact algebraic identities
\[
 \mathsf P_++\mathsf P_-=I,\qquad \mathsf P_\pm^2=\mathsf P_\pm+E,\qquad
 \mathsf P_+\mathsf P_-=\mathsf P_-\mathsf P_+=-E.
\]
In particular, the operators are not projections when $E\ne0$.
For $u_\sigma=\mathsf P_\sigma\psi$, $\sigma\in\{+1,-1\}$, direct application
of $\mathsf P_\sigma$ gives
\begin{align}
 (D_t+\sigma Q)u_\sigma
 ={}&-2\sigma QE\psi-[\mathsf P_\sigma,D_t]\psi
       -[\mathsf P_\sigma,Q](\mathsf P_+-\mathsf P_-)\psi,
 \qquad \psi=u_++u_-.
 \label{eq:HHexactcoupled}
\end{align}
Thus the defect term acts on both components.  Reapplying $\mathsf P_\sigma$
does not eliminate this dependence: it introduces the terms involving
$(D_t+\sigma Q)E\psi$ displayed in (2.12)--(2.13) of~\cite{HH2026}.
On the following page, a sum involving both components is replaced by
the two separate equations
\[
 (D_t\pm Q)(I+\widetilde E)\mathsf P_\pm\psi=0,
\]
where the notation for the defect has been rescaled.  A vanishing sum
does not imply that its two summands vanish.  Here the missing
justification cannot be supplied by orthogonality, since
$\mathsf P_+\mathsf P_-=-E$.

The failure is visible without any pseudodifferential remainder.
Freeze the principal coefficients and suppress the connection terms.
Let
\[
 b(\xi)=\alpha^a e_a^{\ j}\xi_j+M\beta,\qquad
 q_g(\xi)=(M^2+g^{ij}\xi_i\xi_j)^{1/2},\qquad
 q_0(\xi)=(M^2+|\xi|^2)^{1/2},
\]
so that $b(\xi)^2=q_g(\xi)^2I$.  The symbols of the operators above
then satisfy
\[
 p_\pm=\tfrac12\bigl(I\pm b/q_0\bigr),\qquad
 p_\pm^2-p_\pm
   =\frac{q_g^2-q_0^2}{4q_0^2}I
   =\frac{(g^{ij}-\delta^{ij})\xi_i\xi_j}{4(M^2+|\xi|^2)}I.
\]
This is generically of order zero, not of negative order.  Accordingly,
$QE$ is a first-order term.  If $b(\xi)v=q_g(\xi)v$, then the frozen
Dirac solution $\psi=e^{ix\cdot\xi-itq_g(\xi)}v$ has
\[
 p_+\psi=\tfrac12(1+q_g/q_0)\psi.
\]
All coefficient commutators vanish in this frozen calculation.  For
any invertible constant correction $I+\widetilde e$, the proposed plus
equation instead has residual
\[
 (D_t+q_0)(I+\widetilde e)p_+\psi
   =(q_0-q_g)(I+\widetilde e)
       \tfrac12(1+q_g/q_0)\psi,
\]
which is nonzero whenever $q_g(\xi)\ne q_0(\xi)$ and $v\ne0$.
At high frequency the discrepancy is of first order.  An invertible
order-zero normalization of the dependent variable cannot change the
scalar principal dispersion from $q_0$ to $q_g$.

The constant-coefficient calculation is a test of the asserted local
symbolic reduction, not a counterexample to a scattering theorem for
globally asymptotically flat metrics.  The same principal-symbol
obstruction occurs at a point in a smooth compactly supported metric
perturbation.  A valid reduction uses the metric-adapted dispersion
$q_g$ and retains the lower-order couplings created by quantization,
time dependence, and the connection.  Estimating those couplings is a
separate analytical step.  The estimates established here perform that
step under our stated hypotheses; they do not settle the original
theorem under all the hypotheses of~\cite{HH2026}.
\end{remark}

\subsection{The geometric Hamiltonian}
\label{sec:geometry}
Consider $\mathbf g=-dt^2+\gamma_{ij}(t,x)dx^idx^j$ with
$G=\gamma^{-1}=I_3+h$. The frame
$E_0=\partial_t$, $E_a=e_a^{\ j}\partial_j$, $e=G^{1/2}$, is
orthonormal. Its dual coframe has spatial matrix $e^{-1}$.
All these matrices and their inverses are uniformly bounded.
Differentiation of the matrix square root, or its smooth functional
calculus on a fixed compact set of positive matrices, gives
\begin{align*}
 |\partial_x^\gamma(e-I_3)|+
 |\partial_x^\gamma(e^{-1}-I_3)|
 &\lesssim\varepsilon\langle x\rangle^{-\rho-|\gamma|},\\
 |\partial_t\partial_x^\gamma e|+
 |\partial_t\partial_x^\gamma e^{-1}|
 &\lesssim\varepsilon\langle x\rangle^{-\rho-1-|\gamma|}.
\end{align*}
The same spatial estimates hold for $\gamma-I_3$.

For a synchronous metric the possibly nonzero Christoffel coefficients
are combinations of
\[
 \Gamma^0_{ij}=\tfrac12\partial_t\gamma_{ij},\qquad
 \Gamma^k_{0j}=\tfrac12g^{k\ell}\partial_t\gamma_{\ell j},\qquad
 \Gamma^k_{ij}=\tfrac12g^{k\ell}
 (\partial_i\gamma_{j\ell}+\partial_j\gamma_{i\ell}
                                      -\partial_\ell\gamma_{ij}).
\]
Thus each connection coefficient, after any fixed number of spatial
derivatives, has the required decay with one additional spatial power.
In tetrad notation the spin connection is
\[
 \Omega_\mu=\frac14\omega_{\mu AB}\Gamma^A\Gamma^B,
 \qquad
 \omega_{\mu AB}=\mathbf g(E_A,\nabla_{\partial_\mu}E_B),
\]
with constant Clifford matrices $\Gamma^A$ in a fixed convention;
see~\cite{Alcubierre}. Consequently
\[
 |\partial_x^\gamma\Omega_\mu(t,x)|
 \lesssim\varepsilon\langle x\rangle^{-\rho-1-|\gamma|}.
\]
Putting the Dirac equation into Hamiltonian form gives the principal
part $\alpha^ae_a^{\ j}D_j+M\beta$ and zeroth-order combinations of
$\Omega_0$, $e_a^{\ j}\Omega_j$, and fixed matrices. Passing to the
symmetrized expression \eqref{main:hamiltonian} adds only a multiple
of $\alpha^a\partial_je_a^{\ j}$. Hence the resulting $C$ obeys
$\mathfrak c_{\rho+1,J}(C)\lesssim\mathfrak g_{\rho,J+1}(h)$.
This proves the coefficient assertion in
Corollary~\ref{main:geometric}.

In particular, the temporal rotation of the spatial frame is retained
in $C$. There is no time-integrated frame transformation, and no
spatial derivatives of such a transformation have to be controlled.
The geometric Hamiltonian may be symmetric in a density-weighted
space rather than in Lebesgue $L^2$; the analysis of
\eqref{main:hamiltonian} permits precisely this zeroth-order
nonselfadjoint part. We use the original spinor variables throughout,
which also keeps the convolution in \eqref{main:hartree} fixed.

\section{A fixed-energy local-energy estimate}
\label{dirac:scalar-input}

We use the scalar Weyl operator from the preceding reduction,
\[
 A(t)=\Op^w(q_g(t)),\qquad
 q_g(t,x,\xi)=\bigl(M^2+g^{ij}(t,x)\xi_i\xi_j\bigr)^{1/2},
 \qquad H=(M^2-\Delta)^{1/2}.
\]
\begin{lemma}[Explicit scalar perturbation]\label{dirac:scalar-symbol-lemma}
Let $G=I+h$ be real symmetric and uniformly elliptic, with
$|\partial_x^\gamma h|\le\varepsilon
\langle x\rangle^{-\rho-|\gamma|}$ for sufficiently many spatial
derivatives. Then
\begin{equation}\label{dirac:scalar-symbol}
 K_0(t):=A(t)-H=\Op^w(k(t)),\qquad
 k\in\varepsilon SG^{1,-\rho},
\end{equation}
uniformly in time. For sufficiently small $\varepsilon$, $A(t)$ is
selfadjoint on $H^1$, and its graph norm is uniformly equivalent to
the $H^1$ norm.
\end{lemma}
\begin{proof}
Put $h_M(\xi)=(M^2+|\xi|^2)^{1/2}$. The exact pointwise formula
\[
 k(t,x,\xi)=q_g(t,x,\xi)-h_M(\xi)
 =\frac{h^{ij}(t,x)\xi_i\xi_j}{q_g(t,x,\xi)+h_M(\xi)}
\]
gives the assertion directly. Uniform ellipticity and $M>0$ imply
$q_g+h_M\asymp_M\langle\xi\rangle$. Differentiating the quotient
shows that each frequency derivative lowers the frequency order by
one. Every spatial derivative either differentiates the numerator
or differentiates a coefficient in the denominator, and gains the
corresponding inverse spatial power. Thus, for all of the finitely
many derivatives needed subsequently,
\[
 |\partial_x^\gamma\partial_\xi^\delta k|
 \le C_{\gamma,\delta}\varepsilon
       \langle x\rangle^{-\rho-|\gamma|}
       \langle\xi\rangle^{1-|\delta|}.
\]
The reality of $k$ makes $K_0$ symmetric on $H^1$.
The calculus gives $\|K_0H^{-1}\|_{2\to2}\le C\varepsilon$.
The relative bound is less than one for small $\varepsilon$, so the
Kato--Rellich theorem proves selfadjointness on $H^1$. The same
inequality gives the claimed graph-norm equivalence.
\end{proof}

Fix $\rho>2$ and choose
\begin{equation}\label{dirac:weight-choice}
 0<\nu<\min\{1,\rho-2\},\qquad
 a=\frac{3+\nu}{2},\qquad w(x)=\langle x\rangle^{-a}.
\end{equation}
In particular, $3/2<a<2$ and $2a<\rho+1$. Define a smooth radial
function $f_\nu$ by $f_\nu(0)=0$ and
\[
 f_\nu'(r)=\int_0^r(1+z^2)^{-(1+\nu)/2}\,dz,
 \qquad B_\nu=\tfrac12(\nabla f_\nu\cdot D+D\cdot\nabla f_\nu).
\]
The gradient of $f_\nu$ is bounded and belongs to $SG^{0,0}$ as a
spatial coefficient. Write $p=(1+\nu)/2$. The radial and tangential
eigenvalues of its Hessian are $f_\nu''(r)$ and $f_\nu'(r)/r$;
since $f_\nu''$ is decreasing, both are at least
$(1+r^2)^{-p}$. Direct differentiation in dimension three gives
\begin{equation}\label{dirac:multiplier-derivatives}
 \nabla^2f_\nu\ge\langle x\rangle^{-1-\nu}I,\qquad
 -\Delta^2f_\nu
 =2p\bigl(5+(3-2p)|x|^2\bigr)(1+|x|^2)^{-p-2}
 \ge c_\nu w^2.
\end{equation}

\begin{lemma}[Coercivity at all frequencies]\label{dirac:free-coercivity}
For $M>0$ and $w$ as in \eqref{dirac:weight-choice},
\begin{equation}\label{dirac:scalar-positive}
 \langle i[H,B_\nu]v,v\rangle
 \ge c_{M,\nu}\|H^{1/2}(wv)\|_2^2
\end{equation}
for $v\in\mathcal S(\mathbb R^3)$.
\end{lemma}
\begin{proof}
The exact differential commutator is
\[
 q_\nu=i[H^2,B_\nu]
 =2D\cdot\nabla^2f_\nu D-\tfrac12\Delta^2f_\nu.
\]
Since $|\nabla w|\lesssim w$ and
$\langle x\rangle^{-(1+\nu)/2}\ge w$, it follows from
\eqref{dirac:multiplier-derivatives} that
\begin{equation}\label{dirac:squared-coercivity}
 \langle q_\nu z,z\rangle
 \ge 2\|\langle x\rangle^{-(1+\nu)/2}\nabla z\|_2^2
      +c_\nu\|wz\|_2^2
 \ge c_{M,\nu}\|H(wz)\|_2^2.
\end{equation}
Let $R_\lambda=(H^2+\lambda)^{-1}$ and
$L=wH^2w^{-1}-H^2$. More explicitly,
\[
 L=2\nabla\log w\cdot\nabla+
       \Delta\log w-|\nabla\log w|^2.
\]
The operator $LH^{-1}$ is bounded on $L^2$, as is $HR_\lambda$
uniformly for $\lambda\ge0$. If $z=wR_\lambda v$, then
$(H^2+\lambda+L)z=wv$, and hence
\[
 HR_\lambda(wv)=(I+HR_\lambda LH^{-1})Hz.
\]
Therefore $\|HwR_\lambda v\|_2\ge c\|HR_\lambda(wv)\|_2$,
with $c>0$ independent of $\lambda$. No invertibility of the factor
$I+HR_\lambda LH^{-1}$ is used in this inequality.
The resolvent representation of the square-root commutator gives
\[
 \langle i[H,B_\nu]v,v\rangle
 =\frac1\pi\int_0^\infty\lambda^{1/2}
       \langle q_\nu R_\lambda v,R_\lambda v\rangle\,d\lambda.
\]
This identity is first justified on Schwartz functions by truncating
the integral and then passing to the limit as a quadratic form.
Combining it with \eqref{dirac:squared-coercivity} and using the
spectral identity
\[
 \int_0^\infty\frac{\lambda^{1/2}h^2}{(h^2+\lambda)^2}\,d\lambda
 =\frac\pi2h\qquad(h>0)
\]
proves \eqref{dirac:scalar-positive}.
\end{proof}

Fix $s\ge0$. For sufficiently small fixed $\kappa>0$, set
\[
 E_s^\pm=H^{s-1/2}(H\pm\kappa B_\nu)H^{s-1/2}.
\]
The bounded selfadjoint operator $H^{-1/2}B_\nu H^{-1/2}$ gives
\begin{equation}\label{dirac:energy-equivalence}
 c_s\|u\|_{H^s}^2\le\langle E_s^\pm u,u\rangle
                  \le C_s\|u\|_{H^s}^2.
\end{equation}
These energies are independent of time. They have orders $(2s,0)$,
and scalar symbol cancellation in the commutator with
\eqref{dirac:scalar-symbol} yields
\[
 i[K_0,E_s^\pm]\in\varepsilon SG^{2s,-\rho-1}.
\]
Thus
\[
 H^{-s}w^{-1}i[K_0,E_s^\pm]w^{-1}H^{-s}
 \in\varepsilon SG^{0,2a-\rho-1},
\]
whose $L^2$ norm is $O(\varepsilon)$ since $2a<\rho+1$. We obtain
\begin{equation}\label{dirac:perturbative-commutator}
 |\langle i[K_0,E_s^\pm]u,u\rangle|
 \le C_s\varepsilon\|wu\|_{H^s}^2.
\end{equation}

The following norm equivalence makes the coercivity at regularity
$s$ explicit:
\begin{equation}\label{dirac:weighted-bessel}
 \|H^{1/2}wH^{s-1/2}u\|_2\asymp\|H^s(wu)\|_2.
\end{equation}
Indeed the operators
\[
 T_s=H^{1/2}wH^{s-1/2}w^{-1}H^{-s},\qquad
 S_s=H^swH^{1/2-s}w^{-1}H^{-1/2}
\]
belong to $SG^{0,0}$ and are bounded on $L^2$. On Schwartz functions
their products are exactly the identity, so both sides of
\eqref{dirac:weighted-bessel} control each other. The use of massive
Bessel multipliers ensures that this argument includes the zero
frequency region.

For $(D_t+\sigma A(t))u=0$, the free coercivity, the weighted
equivalence, and \eqref{dirac:perturbative-commutator} show that, for
$\varepsilon$ sufficiently small,
\begin{equation}\label{dirac:energy-monotonicity}
 \frac d{dt}\langle E_s^{-\sigma}u,u\rangle
 \le-c_s\|wu\|_{H^s}^2,\qquad
 \frac d{dt}\langle E_s^{\sigma}u,u\rangle
 \ge c_s\|wu\|_{H^s}^2.
\end{equation}
The first inequality is integrated forward from the initial time,
and the second backward. No time derivative of $A$ occurs in this
calculation.

\section{Forced energy and the full Dirac evolution}
\label{dirac:energy}

\begin{lemma}[Forced scalar estimate]\label{dirac:forced-scalar}
Fix $s\ge0$, $\rho>2$, and \eqref{dirac:weight-choice}. Under the
spatial coefficient bounds above, with sufficiently small
$\varepsilon$, solutions of
\[
 (D_t+\sigma A(t))u=F_1+F_2,
 \qquad \sigma\in\{+1,-1\},
\]
on an interval $I$ with initial endpoint $t_0$ satisfy
\begin{equation}\label{dirac:forced-LE}
 \|u\|_{L^\infty(I;H^s)}+\|wu\|_{L^2(I;H^s)}
 \le C_s\bigl(\|u(t_0)\|_{H^s}
       +\|F_1\|_{L^1(I;H^s)}
       +\|w^{-1}F_2\|_{L^2(I;H^s)}\bigr).
\end{equation}
The estimate holds for either orientation of $I$, with a constant
independent of its length and position.
\end{lemma}
\begin{proof}
On a forward interval use $E_s^{-\sigma}$ in
\eqref{dirac:energy-monotonicity}. The source contribution is at most
$2|\langle E_s^{-\sigma}u,F_1+F_2\rangle|$. The first pairing is
bounded by $C_s\|u\|_{H^s}\|F_1\|_{H^s}$. For the second, put
$z=wu$ and $f=w^{-1}F_2$. Then
\[
 \langle E_s^{-\sigma}u,F_2\rangle
 =\langle wE_s^{-\sigma}w^{-1}z,f\rangle.
\]
The conjugated operator has orders $(2s,0)$, so
$H^{-s}wE_s^{-\sigma}w^{-1}H^{-s}$ is bounded on $L^2$.
Consequently the second pairing is bounded by
$C_s\|wu\|_{H^s}\|w^{-1}F_2\|_{H^s}$.
Integrate the energy inequality, take the supremum in its upper
endpoint, and use Cauchy--Schwarz in time. With
$X=\|u\|_{L^\infty H^s}$ and $Y=\|wu\|_{L^2H^s}$, this gives
\[
 X^2+Y^2\le C\|u(t_0)\|_{H^s}^2
       +CX\|F_1\|_{L^1H^s}+CY\|w^{-1}F_2\|_{L^2H^s}.
\]
Young's inequality proves \eqref{dirac:forced-LE}. On a backward
interval use $E_s^\sigma$ and the second inequality in
\eqref{dirac:energy-monotonicity}. The estimate on an unbounded
interval follows by exhaustion from finite intervals.
\end{proof}

\begin{proposition}[Full Dirac local energy]\label{dirac:full-LE}
Assume \eqref{main:coeff} with $\rho>2$, and choose $w$ as in
\eqref{dirac:weight-choice}. For every fixed $s\ge0$, sufficiently
many controlled coefficient derivatives and sufficiently small
$\varepsilon$ give
\begin{equation}\label{dirac:homogeneous-LE}
 \sup_{t,r}\|U_D(t,r)\|_{H^s\to H^s}\le C_s,
 \qquad
 \|wU_D(\cdot,r)f\|_{L_t^2H_x^s}\le C_s\|f\|_{H^s}.
\end{equation}
The forced estimate \eqref{dirac:forced-LE} also holds with
$(D_t+\sigma A)u$ replaced by $(D_t+B)\psi$.
\end{proposition}
\begin{proof}
Set $u_\sigma=P_\sigma\psi$ and use the exact coupled equations
\eqref{dirac:coupled}. Each remainder satisfies
$R_\sigma\in\varepsilon SG^{0,-\rho-1}$; hence
\[
 w^{-1}R_\sigma w^{-1}
 \in\varepsilon SG^{0,2a-\rho-1}\subset\varepsilon\Psi^0
\]
and, after conjugation by $H^s$,
\begin{equation}\label{dirac:coupling-absorption}
 \|w^{-1}R_\sigma v\|_{L^2H^s}
 \le C_s\varepsilon\|wv\|_{L^2H^s}.
\end{equation}
The operators $P_\sigma$ and $w^{-1}P_\sigma w$ have orders
$(0,0)$ and are uniformly bounded on $H^s$. Apply
Lemma~\ref{dirac:forced-scalar} to both component equations and use
the exact identity $\psi=u_++u_-$. The result is
\[
 \sum_\sigma\bigl(\|u_\sigma\|_{L^\infty H^s}
                   +\|wu_\sigma\|_{L^2H^s}\bigr)
 \le C_s\bigl(\|\psi(t_0)\|_{H^s}
       +\|F_1\|_{L^1H^s}+\|w^{-1}F_2\|_{L^2H^s}\bigr)
       +C_s\varepsilon\sum_\sigma\|wu_\sigma\|_{L^2H^s}.
\]
The last term is absorbed for small $\varepsilon$. Summation of
the components proves the forced estimate and
\eqref{dirac:homogeneous-LE}. In particular, the nonselfadjoint
zeroth-order part of the Hamiltonian is included in the
remainder and is controlled by \eqref{dirac:coupling-absorption}.

We finish with the evolution and approximation details. For smooth
coefficients the Weyl principal part of $B$ is symmetric with
Hermitian order-one symbol. On each bounded time interval, its
ordinary $H^r$ energy estimate follows from
$[H^r,\Op^w(b_M)]H^{-r}\in\Psi^0$ after conjugation of the
equation; the order-zero coefficient $C(t,x)$ is bounded on $H^r$.
Frequency regularization and this estimate construct the evolution
on smooth data. One first obtains convergence with one additional
Sobolev derivative, then uses the uniform finite-interval estimate
and density to obtain a strongly continuous evolution on $H^r$.
The cocycle and inverse identities follow from uniqueness. The same
construction applies to the scalar $A(t)$. Only continuity in time
of the lower-order coefficients is used: differentiating the spin
connection, which would introduce a second time derivative of the
metric, is unnecessary.

The preceding estimates can first be justified on smooth solutions
by spatial regularization. For finitely regular metrics, smooth
space-time approximations preserve \eqref{main:coeff} up to a fixed
constant. They converge locally in time in the required unweighted
coefficient norms after reserving a finite derivative margin. The
corresponding generators converge strongly from $H^{r+1}$ to $H^r$;
finite-interval variation of constants and density give convergence
of the propagators on $H^r$. Weighted local energy passes on bounded
intervals by lower semicontinuity and on the entire line by monotone
convergence. Choose the finite derivative count with a margin for
$r=s+1$ in this construction. The resulting bounds are uniform in
the length and position of the time interval.
\end{proof}

\section{Fourier localizations and free Dirac estimates}
\label{dirac:fourier}
Let $\D=\{1,2,4,\ldots\}$. Choose smooth inhomogeneous Fourier cutoffs
$P_N=p_N(D)$ forming a partition of unity, with $P_1$ the low-frequency
ball and $P_N$ supported where $|\xi|\asymp N$ for $N\ge2$.
For every $r\in\R$,
\begin{equation}\label{td:LP}
 \|f\|_{H^r}^2\asymp\sum_{N\in\D}N^{2r}\|P_Nf\|_2^2.
\end{equation}
For $\mu\in\D$, use a smooth bounded-overlap lattice partition
$\{P_q:q\in\mathcal Q_\mu\}$ of frequency space into cubes of side
comparable to $\mu$. Its convolution kernels have uniformly bounded
$L^1$ norms, and $\sum_q\|P_qf\|_2^2\lesssim\|f\|_2^2$.
Smooth cutoffs are never treated as idempotent projections; enlarged
cutoffs equal to one on their supports are inserted when needed.
Set
\begin{align}
 \Gamma_{\mu,N}(u)(t)^2
 &=\sum_{q\in\mathcal Q_\mu}\|P_qP_Nu(t)\|_\infty^2,\label{td:Gamma}\\
 \mathfrak b_N(v_N)(t)^2
 &=\sum_{\mu\in\D}\mu^{-2}\sum_{q\in\mathcal Q_\mu}
       \|P_qv_N(t)\|_\infty^2,\qquad
 B_N(u)=\mathfrak b_N(P_Nu).\label{td:B}
\end{align}
The notation $\mathfrak b_N$ acts directly on the supplied function
$v_N$ supported in the $N$ frequency block; it does not apply $P_N$ again.
For any fixed exponent $\alpha$ with $1/2<\alpha<s$, write
\begin{equation}\label{td:alpha}
 \epsilon=\alpha-\tfrac12>0,\qquad \delta=s-\alpha>0.
\end{equation}
For an interval $I$, define
\begin{align}
 Q_s(u;I)^2
 &=\sum_{N\in\D}N^{2(s-\alpha)}\|B_N(u)\|_{L^2(I)}^2,\label{td:Q}\\
 \|u\|_{W^s(I)}&=\|u\|_{L^\infty(I;H^s)}+Q_s(u;I).\label{td:W}
\end{align}
The square $Q_s(u;I)^2$ is the integral over $I$ of a nonnegative
density. In particular it is absolutely continuous under time
restriction whenever finite.

Bernstein and cube almost orthogonality imply
\[
 \Gamma_{\mu,N}(u)^2\lesssim\min(\mu,N)^3\|P_Nu\|_2^2.
\]
Since
\[
 \sum_{\mu\le N}\mu^{-2}\mu^3+
 \sum_{\mu>N}\mu^{-2}N^3\lesssim N,
\]
we obtain the useful purely spatial estimates
\begin{equation}\label{td:Bernstein}
 \mathfrak b_N(v_N)\lesssim N^{1/2}\|v_N\|_2,
 \qquad Q_s(u;I)\lesssim\|u\|_{L^2(I;H^{s-\epsilon})}.
\end{equation}
These inequalities do not use an equation or a spatial support condition.

Let $B_0=\alpha\cdot D+M\beta$, $S_D(t)=e^{-itB_0}$ and
\[
 \Pi_\sigma^0(\xi)=\tfrac12\left(I+
  \sigma\frac{\alpha\cdot\xi+M\beta}{\sqrt{M^2+|\xi|^2}}\right).
\]
Then $S_D(t)=\sum_\sigma e^{-i\sigma tH}\Pi_\sigma^0(D)$.
The Fourier projectors preserve every frequency block. Their symbol
derivatives through order two are bounded, since $M>0$, and the exact
moment formula
\[
 x^\gamma\Pi_\sigma^0(D)F
 =\sum_{\delta\le\gamma}\binom\gamma\delta i^{|\delta|}
   (\partial_\xi^\delta\Pi_\sigma^0)(D)x^{\gamma-\delta}F
\]
shows
$\|\langle x\rangle^2\Pi_\sigma^0F\|_2
\lesssim_M\|\langle x\rangle^2F\|_2$.

\begin{lemma}[Retarded truncation with two integrable delay regimes]
\label{paper:whitney}
Let $E_N$ be either the normed space with norm $\mathfrak b_N$ on the
$N$ frequency block, or that block with its $L^\infty_x$ norm. Let
$S(t)=e^{-i\sigma t\sqrt{M^2-\Delta}}$, $M>0$, and suppose
\[
 \|S(t)f_N\|_{L_t^2E_N}\le A_N\|f_N\|_2,\qquad
 \left\|\int_{\mathbb R}S(-r)F_N(r)\,dr\right\|_2
 \le C\|\langle x\rangle^2F_N\|_{L^2_{t,x}}.
\]
For the cube norm take $A_N=C_\eta N^{1/2+\eta}$, and for the
annular norm take $A_N=C_\eta N^{1+\eta}$. Then, for any larger
positive frequency loss, the same respective estimate holds for
$\int_{r<t}S(t-r)F_N(r)\,dr$.
\end{lemma}
\begin{proof}
Let $\mathcal W_j$ consist of dyadic rectangles $I\times J$ with
$|I|=|J|=2^j$, $I$ to the right of $J$, with $I,J$ nonadjacent but
their parents adjacent. These rectangles partition $\{r<t\}$ up
to a null set; on each rectangle $2^j\lesssim t-r\lesssim 2^j$.
Each interval has a uniformly bounded number of partners. Define
\[
 T_NF(t)=S(t)\int_{\mathbb R}S(-r)F(r)\,dr,
 \quad T_{N,j}F=\sum_{I\times J\in\mathcal W_j}
                       \mathbf1_I T_N(\mathbf1_JF).
\]
The hypotheses give $\|T_NF_N\|_{L_t^2E_N}
\le CA_N\|\langle x\rangle^2F_N\|_2$. Since the output intervals
at a fixed scale are disjoint and partner multiplicities are bounded,
\begin{align*}
 \|T_{N,j}F_N\|_{L_t^2E_N}^2
 &\lesssim\sum_{I\times J\in\mathcal W_j}
       \|\mathbf1_I T_N(\mathbf1_JF_N)\|_{L_t^2E_N}^2\\
 &\lesssim A_N^2\sum_{I\times J\in\mathcal W_j}
       \|\mathbf1_J\langle x\rangle^2F_N\|_2^2
 \lesssim A_N^2\|\langle x\rangle^2F_N\|_2^2.
\end{align*}
This estimate is independent of $j$.

For $E_N$ equal to the cube norm, sum these bounds over scales
$1\lesssim2^j\lesssim N^8$. There are $O(1+\log N)$ such scales.
On the union of smaller rectangles, Minkowski, Bernstein, and the
positive delay kernel $\mathbf1_{0<\tau\lesssim1}$ give the bound
$CN^{1/2}\|F_N\|_{L^2_{t,x}}$. On the union of larger rectangles,
use the annular dispersive kernel and count the cubes meeting the
$N$ block:
\[
 \mathfrak b_N(S(\tau)F_N)
 \lesssim N^{5/2}|\tau|^{-3/2}
 \left(\sum_{\mu\ge1}\mu^{-2}(1+N/\mu)^3\right)^{1/2}\|F_N\|_1
 \lesssim N^4|\tau|^{-3/2}\|\langle x\rangle^2F_N\|_2.
\]
The last positive time kernel has bounded integral for
$\tau\gtrsim N^8$. The actual Whitney masks are bounded pointwise
by the indicated delay masks before applying these scalar convolution
estimates; no translation invariance of the Whitney union is required.

For $E_N=L^\infty_x$, use a dyadic lower cutoff comparable to
$N^{-1/2}$ and an upper cutoff comparable to $N^6$. The near-delay
bound is $N^{3/2}N^{-1/2}=N$. The far kernel is
$N^{5/2}|\tau|^{-3/2}\mathbf1_{\tau\gtrsim N^6}$ and has integral
$O(N^{-1/2})$. The intermediate range again has $O(1+\log N)$
scales. Absorb the logarithm by an arbitrarily small additional
frequency power. Smooth compact-time sources suffice initially;
the estimate defines the retarded operator by completion. Reversing
time proves the advanced estimate.
\end{proof}

\begin{lemma}\label{dirac:free}
For each $\eta>0$ and each frequency block $N\ge1$,
\begin{align}
 \|\mathfrak b_N(S_Df_N)\|_{L_t^2}
 &\lesssim_\eta N^{1/2+\eta}\|f_N\|_2,\label{dirac:freecube}\\
 \|S_Df_N\|_{L_t^2L_x^\infty}
 &\lesssim_\eta N^{1+\eta}\|f_N\|_2,\label{dirac:freeann}\\
 \left\|\mathfrak b_N\left(\int_{r<t}S_D(t-r)F_N(r)\,dr\right)
       \right\|_{L_t^2}
 &\lesssim_\eta N^{1/2+\eta}\|\langle x\rangle^2F_N\|_{L^2_{t,x}},
       \label{dirac:retcube}\\
 \left\|\int_{r<t}S_D(t-r)F_N(r)\,dr\right\|_{L_t^2L_x^\infty}
 &\lesssim_\eta N^{1+\eta}\|\langle x\rangle^2F_N\|_{L^2_{t,x}}.
       \label{dirac:retann}
\end{align}
The retarded estimates also hold in the reverse time orientation.
\end{lemma}
\begin{proof}
It suffices to prove the scalar half-wave estimates and use the
projector bounds above. The free cube kernel estimate in
\cite[Lemma 3.6 and Remark 3.7]{HT} has the form
\[
 |K_{N,\mu}(t,x)|\lesssim_M
 \min\{\mu^3,\mu N|t|^{-1},N^{5/2}|t|^{-3/2}\},\quad1\le\mu\le N.
\]
Its time integral is bounded by $C\mu N(1+\log N)$, splitting at
$N/\mu^2$ and $N^3/\mu^2$. The $TT^*$ bound, cube overlap and
summation with weight $\mu^{-2}$ give
$\|\mathfrak b_N(Sf_N)\|_{L^2}^2\lesssim N(1+\log N)\|f_N\|_2^2$.
For $\mu>N$, use bounded overlap with scale-$N$ cubes. Taking
$\mu\asymp N$ also proves the full-annulus estimate. The low block
uses the massive integrable $|t|^{-3/2}$ tail.

The unperturbed scalar energy calculation implies
\[
 \|\langle x\rangle^{-2}S(t)f\|_{L^2_{t,x}}\lesssim\|f\|_2,
 \quad
 \left\|\int S(-r)F(r)\,dr\right\|_2
 \lesssim\|\langle x\rangle^2F\|_{L^2_{t,x}}.
\]
Apply Lemma~\ref{paper:whitney} with an initial homogeneous loss
$0<\eta_0<\eta$. Its logarithmic factor is absorbed by
$N^{\eta-\eta_0}$. The flat Dirac projectors preserve frequency
blocks and the weighted source norm, so summing the two scalar
estimates proves the retarded assertions.
\end{proof}

\section{Frequency-dependent exterior gluing}
\label{sec:gluing}
Write $B=B_0+K$, where $K$ includes the principal perturbation
and all zeroth-order terms. The following summation lemma applies to matrix symbols componentwise.
\subsection{A square-sum lemma for varying operator families}

\begin{lemma}\label{declin:family}
Let $T_N$, $N\in\{1,2,4,\ldots\}$, be a family of spatial
pseudodifferential operators whose symbols are uniformly bounded in
$S^m_{1,0}$. For any fixed real $r$,
\begin{equation}\label{declin:familybound}
 \sum_NN^{2(r-m)}\|P_NT_Nz\|_2^2
 \le C_{r,m}\|z\|_{H^r}^2.
\end{equation}
The same assertion holds with $P_N$ replaced by any uniformly smooth
annular Fourier multiplier at scale $N$. The constant requires only
finitely many uniform symbol seminorms, depending on $r,m$.
\end{lemma}

\begin{proof}
We give the off-diagonal estimate because a uniform operator bound
for each $T_N$ alone would not justify the square sum.
For each integer $J>0$, finite-order symbol composition and the
$L^2$ boundedness theorem give, uniformly in $N,M$,
\[
 \|T_NP_M\|_{2\to2}\lesssim M^m,
 \qquad
 \|\langle D\rangle^J T_NP_M\|_{2\to2}
 \lesssim M^{m+J}.
\]
In the second bound, first compose with $\langle D\rangle^J$ and
then use the input multiplier $P_M$. The resulting symbol is supported
in the $M$ block and its finitely many normalized seminorms are
bounded by $CM^{m+J}$. If $N\ge4M$, this proves
\[
 \|P_NT_NP_M\|_{2\to2}\lesssim N^{-J}M^{m+J}.
\]
For $M\ge4N$, apply the same argument to $T_N^*$, whose symbols
have uniformly bounded seminorms in $S^m_{1,0}$, and then take the
adjoint. On comparable blocks use the first bound. Increasing $J$
as necessary, we have for every prescribed $L>0$
\begin{equation}\label{declin:offdiag}
 \|P_NT_NP_M\|_{2\to2}
 \le C_L M^m2^{-L|\log_2N-\log_2M|}.
\end{equation}
Enlarged smooth cutoffs can be inserted throughout, so no cutoff is
being treated as an idempotent projection.

Expand $z=\sum_MP_Mz$ and multiply \eqref{declin:offdiag} by
$N^{r-m}$. With $z_M=M^r\|P_Mz\|_2$, the resulting sequence is
bounded by convolution of $z_M$ with
$C2^{-(L-|r-m|)|\log_2N-\log_2M|}$.
Choose $L>|r-m|+1$ and use the discrete Young inequality and
Littlewood--Paley equivalence. This proves
\eqref{declin:familybound}. All compositions and adjoints above use
only finitely many spatial derivatives; choosing the corresponding
finite $J,L$ fixes the required seminorm count.
\end{proof}

\subsection{Cutoff powers and the Hartree resolution estimate}
Fix $\chi\in C_c^\infty$ equal to one on $|x|\le1$ and zero on
$|x|\ge2$, and put $\chi_R(x)=\chi(x/R)$. For $\rho>a+2$,
the following normalized symbol families are uniform in $t,R\ge1$:
\begin{align*}
 R^{-a}\chi_Rw^{-1}&\in\Psi^0,\\
 R^{-(a+1)}x^\gamma[B_0,\chi_R]w^{-1}&\in\Psi^0,
                    &&|\gamma|\le2,\\
 R^{\rho-a-2}x^\gamma(1-\chi_R)Kw^{-1}&\in\varepsilon\Psi^1,
                    &&|\gamma|\le2.
\end{align*}
Indeed $[B_0,\chi_R]=-i\alpha\cdot\nabla\chi_R$ is supported on
$|x|\asymp R$ and has size $R^{-1}$. The factor $x^\gamma w^{-1}$
costs at most $R^{a+2}$ there. For the tail, the spatial symbol
order is $a+2-\rho<0$ outside $|x|\ge R$.
Commuting the moments through $P_N$ uses the same exact polynomial
identity as above; differentiated cutoff symbols gain $N^{-1}$.
Thus Lemma~\ref{declin:family}, with $R_N=N^\theta$, gives
\begin{align}
 \sum_NN^{2(s-a\theta)}\|P_N\chi_{R_N}u\|_{L^2_{t,x}}^2
 &\lesssim\|wu\|_{L_t^2H^s}^2,\label{dirac:inside-square}\\
 \sum_NN^{2(s-(a+1)\theta)}
 \|\langle x\rangle^2P_N[B_0,\chi_{R_N}]u\|_{L^2_{t,x}}^2
 &\lesssim\|wu\|_{L_t^2H^s}^2,\label{dirac:boundary-square}\\
 \sum_NN^{2(s-1+\theta(\rho-a-2))}
 \|\langle x\rangle^2P_N(1-\chi_{R_N})Ku\|_{L^2_{t,x}}^2
 &\lesssim\varepsilon^2\|wu\|_{L_t^2H^s}^2.\label{dirac:tail-square}
\end{align}
For $u=U_D(t,r)f$ and $v_N=P_N((1-\chi_{R_N})u)$, the exact equation is
\[
 (i\partial_t-B_0)v_N
 =P_N\big((1-\chi_{R_N})Ku+[B_0,\chi_{R_N}]u\big).
\]
No Fourier projection has been commuted through the variable flow.

Choose
\begin{equation}\label{dirac:hartree-parameters}
 \frac{a+1}{\rho-1}<\ell<\min\{1,s-\tfrac12\},\quad
 \theta=\frac\ell{a+1},\quad
 0<\eta<s-\tfrac12-\ell,\quad \alpha=\tfrac12+\ell+\eta.
\end{equation}
These parameters are available in the Hartree range of
Theorem~\ref{main:hartree-theorem}, taking $a>3/2$ sufficiently close to $3/2$.
The interior estimate \eqref{td:Bernstein} and
\eqref{dirac:inside-square} cost $N^{1/2+a\theta}$.
The exterior source losses in \eqref{dirac:boundary-square} and
\eqref{dirac:tail-square} are
\[
 (a+1)\theta=\ell,\qquad
 1-\theta(\rho-a-2)<\ell.
\]
Lemma~\ref{dirac:free} therefore gives
\begin{equation}\label{dirac:hartree-flow}
 \sup_r\|U_D(\cdot,r)f\|_{W^s}\lesssim_s\|f\|_{H^s}.
\end{equation}
For clarity, multiplying the cube output by $N^{s-\alpha}$ leaves
the boundary source weight $N^{s-\ell}$, exactly the exponent
allowed above, and leaves the interior weight
$N^{s-\ell-\eta}$, which is at most $N^{s-a\theta}$.
The free initial exterior term is handled by
Lemma~\ref{declin:family} for $1-\chi_{R_N}$ at index $s-\ell$.
Use Proposition~\ref{dirac:full-LE} for the remaining weighted
Sobolev norms. This proves \eqref{dirac:hartree-flow} with the actual
full Dirac energy input, including all connection terms.

\section{A local Strichartz improvement for the cubic Dirac equation}
\label{dirac:local-cubic}

Throughout this section $B(t)=B_0+K(t)$ is the full Dirac Hamiltonian,
$B_0=\alpha\cdot D+M\beta$, and $U_D(t,r)$ is its exact evolution.
All connection terms are included in $K$.  We use the uniform
energy and weighted local-energy estimate of Proposition~\ref{dirac:full-LE},
\begin{equation}\label{dirac:cubic-LE}
 \sup_t\|U_D(t,r)f\|_{H^s}
 +\|w_aU_D(\cdot,r)f\|_{L^2_tH^s_x}
 \le C_s\|f\|_{H^s},\qquad w_a=\langle x\rangle^{-a},
\end{equation}
where $3/2<a<2$, together with the stated small polynomial coefficient
bounds.  In particular, the principal coefficients and their first time
derivatives have the required spatial bounds; the zeroth-order matrix
$C$ needs only its uniform spatial bounds and time continuity.
No second time derivative is used in this section.

\begin{lemma}[Uniform local estimate]\label{dirac:local-stz}
For every $\eta>0$ there is $T_0>0$ such that, on every interval
$I$ of length at most $T_0$ and for every $t_0\in I$,
\begin{equation}\label{dirac:local-stz-bound}
 \|v\|_{L^2(I;L^\infty)}
 \le C_\eta\left(\|v(t_0)\|_{H^{1+\eta}}
       +\|(D_t+B)v\|_{L^1(I;H^{1+\eta})}\right).
\end{equation}
The constants are independent of the position of $I$.
\end{lemma}

\begin{proof}
We reduce to a local wave estimate and retain the order-zero coupling.
The local wave parametrix theorem of Smith~\cite{Smith1998}, in particular
the Pecher estimates in its Section~6, p.~827, gives in dimension three
\[
 \|z\|_{L^p(I;L^q)}\lesssim
 \|z(t_0)\|_{H^1}+\|\partial_tz(t_0)\|_2
 +\|(\partial_t^2-g^{ij}\partial_i\partial_j)z\|_{L^1(I;L^2)},
 \qquad q>6,\quad p=\frac{2q}{q-6}.
\]
This theorem permits coefficients Lipschitz in time, with two bounded
spatial derivatives (the precise assumptions are on p.~798).
The interval length and constant can be chosen
uniformly under our coefficient bounds.  Applying the estimate to
$\langle D\rangle^\eta z$ is legitimate with our additional spatial
regularity: the commutator with the spatial second-order operator has
order $1+\eta$ and is controlled by
$\|z\|_{L^1(I;H^{1+\eta})}$.

Let $A=\operatorname{Op}^w(q_g)$, with
$q_g=(M^2+g^{ij}\xi_i\xi_j)^{1/2}$, and first consider a homogeneous
scalar solution $(D_t+\sigma A)z=0$, $\sigma=\pm1$.
It satisfies
\[
 \partial_t^2z+A^2z=-i\sigma(\partial_tA)z.
\]
Here $\partial_tA$ has order one, and
$A^2+g^{ij}\partial_i\partial_j$ has order at most one.
Consequently the source in the second-order wave estimate is bounded in
$L^1H^\eta$ by $C|I|\sup_I\|z\|_{H^{1+\eta}}$.
The local Sobolev energy estimate controls this supremum by the initial
$H^{1+\eta}$ norm.  Choose a finite $q>\max\{6,3/\eta\}$.
The Sobolev embedding $W^{\eta,q}\hookrightarrow L^\infty$, followed by
$L^p(I)\hookrightarrow L^2(I)$, proves the scalar homogeneous version of
\eqref{dirac:local-stz-bound}.  Scalar Duhamel and Minkowski give the
inhomogeneous estimate with an $L^1H^{1+\eta}$ source; in particular,
no time derivative of that source is taken.

Apply this scalar estimate to the corrected components
$P_\sigma v$.  Their exact equations are
\[
 (D_t+\sigma A)P_\sigma v
 =P_\sigma (D_t+B)v+R_\sigma v,
 \qquad R_\sigma\in\Psi^0.
\]
The operators $P_\sigma$ and $R_\sigma$ are uniformly bounded on
$H^{1+\eta}$.  The ordinary finite-interval Sobolev energy estimate for
$B$ bounds $\sup_I\|v\|_{H^{1+\eta}}$ by the initial norm and the
$L^1H^{1+\eta}$ source.  Summing the two estimates and using the exact
identity $v=P_+v+P_-v$ proves the lemma.  We have not used boundedness of
a pseudodifferential operator on $L^\infty$.
\end{proof}

\begin{lemma}[Interior estimate with varying radii]
\label{dirac:cubic-interior}
Let $R_N=N^\theta$, $\theta>0$, and let $u$ solve
$(D_t+B)u=0$.  For $r=1+\eta$, every real $s$, and a fixed smooth
cutoff $\chi_R(x)=\chi(x/R)$, one has
\begin{equation}\label{dirac:cubic-interior-bound}
 \sum_{N\in\mathcal D}N^{2(s-r-a\theta)}
 \|P_N(\chi_{R_N}u)\|_{L^2_tL^\infty_x}^2
 \le C\|w_au\|_{L^2_tH^s_x}^2.
\end{equation}
Only finitely many uniform coefficient derivatives, depending on $s,r$,
are needed.
\end{lemma}

\begin{proof}
Choose intervals $I_j$ covering the line with bounded overlap, and smooth
functions $\zeta_j$ supported in slightly larger intervals $I_j^*$ of
length at most $T_0$, equal to one on $I_j$, and with uniformly bounded
first derivatives.  Put
\[
 z_{j,N}(t)=\zeta_j(t)P_N(\chi_{R_N}u(t)).
\]
It has zero initial data at the left endpoint of $I_j^*$, and satisfies
\[
 (D_t+B)z_{j,N}
 =-i\zeta_j'P_N\chi_{R_N}u
   +\zeta_j[B,P_N\chi_{R_N}]u.
\]
Lemma~\ref{dirac:local-stz}, Cauchy--Schwarz on $I_j^*$, and bounded
overlap imply
\begin{align}
 \|P_N\chi_{R_N}u\|_{L^2_tL^\infty_x}^2
 \lesssim\int_{\mathbb R}\bigl(
  \|P_N\chi_{R_N}u(t)\|_{H^r}^2
 +\|[B,P_N\chi_{R_N}]u(t)\|_{H^r}^2\bigr)\,dt.
 \label{dirac:windowed-interior}
\end{align}

We justify the required frequency square summation, rather than summing
individual operator norms.  Write $z=w_au$ and set
\[
 T_N^{(0)}=R_N^{-a}P_N\chi_{R_N}w_a^{-1},\qquad
 T_N^{(1)}=R_N^{-a}[B,P_N\chi_{R_N}]w_a^{-1}.
\]
For either family and every prescribed $L>0$,
\begin{equation}\label{dirac:interior-matrix}
 \|T_N^{(i)}P_M\|_{L^2\to H^r}
 \le C_L N^r 2^{-L|\log_2N-\log_2M|},\qquad i=0,1.
\end{equation}
For $i=0$, the output is in the $N$ annulus, and the multiplier
$R_N^{-a}\chi_{R_N}w_a^{-1}$ has uniform bounded derivatives of every
required finite order.  The ordinary high-to-low and low-to-high
pseudodifferential estimates prove \eqref{dirac:interior-matrix}.
For $i=1$, use
\[
 [B,P_N\chi_R]=[B,P_N]\chi_R+P_N[B,\chi_R].
\]
The second term again has output frequency $N$, and
$[B,\chi_R]$ is multiplication by a bounded matrix of size $O(R^{-1})$
on the cutoff shell.  The normalized coefficient
$R^{-a}[B,\chi_R]w_a^{-1}$ has uniformly bounded derivatives.
For the first term, the kernel of $[b^j(x),P_N]D_j$ contains
$b^j(x)-b^j(y)$.  A first spatial derivative of $b^j$ gains $N^{-1}$,
which cancels the derivative $D_j$ when the input frequency is
comparable to $N$.  The zeroth-order part of $B$ gives a still lower
order commutator.  Here is an explicit verification that its output
need not be localized to the $N$ annulus.  Let $Q_H$ be an output
Littlewood--Paley projector, and put $S_N=[B,P_N]$.  The block
$Q_HS_NP_K$ vanishes unless $H\simeq N$ or $K\simeq N$, since
$B P_N-P_N B$ contains $P_N$ on one side or the other.
For $H,K\simeq N$ the preceding commutator bound gives order zero.
If $K\simeq N$ and $H\gg N$, a coefficient piece at frequency
$\simeq H$ is required; its $L^\infty$ norm is at most $C_JH^{-J}$,
while the derivative costs $O(N)$.  If $H\simeq N$ and $K\gg N$,
the corresponding bound is $C_JK^{1-J}$.  The cases with one
frequency much below $N$ are bounded similarly using coefficient
frequency $\simeq N$.  Thus, after increasing the finite number of
coefficient derivatives, for any prescribed $A>0$,
\[
 \|Q_HS_NP_K\|_{2\to2}
 \le C_A2^{-A|\log_2(H/N)|}2^{-A|\log_2(K/N)|}.
\]
The multiplier $a_N=R_N^{-a}\chi_{R_N}w_a^{-1}$ satisfies
$\|P_Ka_NP_M\|_{2\to2}\le
C_A2^{-A|\log_2(K/M)|}$ uniformly in $N$.
Convolving these matrices in $K$, and then summing the output blocks
with the $H^r$ weight, proves \eqref{dirac:interior-matrix}; choose
$A$ larger than $L+|r|+2$.  This also makes explicit that the finite
number of derivatives required depends on both $L$ and $r$.

Expand $z=\sum_MP_Mz$, multiply \eqref{dirac:interior-matrix} by
$N^{s-r}$, and use discrete Young's inequality, with
$L>|s|+2$.  It follows that
\[
 \sum_NN^{2(s-r)}\|T_N^{(i)}z\|_{H^r}^2
 \lesssim\|z\|_{H^s}^2.
\]
Since $R_N^a=N^{a\theta}$, applying this estimate pointwise in time
to \eqref{dirac:windowed-interior} proves
\eqref{dirac:cubic-interior-bound}.
\end{proof}

\begin{proposition}[Global cubic resolution estimate]
\label{dirac:cubic-global}
Suppose $\rho>a+2$ and choose $\ell,\theta,\eta>0$ such that
\begin{equation}\label{dirac:cubic-parameters}
 \frac{a+1}{\rho-1}<\ell<\min\{1,s-1\},\qquad
 \theta=\frac\ell{a+1},\qquad 0<\eta<s-1-\ell.
\end{equation}
Then
\begin{equation}\label{dirac:cubic-global-bound}
 \|U_D(\cdot,r)f\|_{L^\infty_tH^s_x\cap L^2_tL^\infty_x}
 \le C_s\|f\|_{H^s},
\end{equation}
uniformly in the initial time $r$.  The corresponding Duhamel operator
maps $L^1_tH^s_x$ to the same space in both time orientations.
\end{proposition}

\begin{proof}
Write $u=U_D(\cdot,r)f$ and split each output frequency as
\[
 P_Nu=P_N\chi_{R_N}u+v_N,\qquad
 v_N=P_N((1-\chi_{R_N})u).
\]
The interior estimate has exponent $s-1-\eta-a\theta>0$ by
\eqref{dirac:cubic-parameters}.  Therefore Cauchy--Schwarz in $N$
and Lemma~\ref{dirac:cubic-interior} give the required
$L^2_tL^\infty_x$ bound for the sum of the interior pieces.

The exterior equation uses the flat Dirac operator exactly:
\[
 (D_t+B_0)v_N=F_N,
\]
where the exact source is
\[
 F_N=-P_N\bigl((1-\chi_{R_N})K u+[B_0,\chi_{R_N}]u\bigr).
\]
The varying-family weighted estimates used for the Hartree argument
apply also to this differential Dirac perturbation.  They give
\begin{align*}
 \sum_NN^{2(s-(a+1)\theta)}
  \|\langle x\rangle^2P_N[B_0,\chi_{R_N}]u\|_{L^2_{t,x}}^2
 &\lesssim\|w_au\|_{L^2_tH^s_x}^2,\\
 \sum_NN^{2(s-1+\theta(\rho-a-2))}
  \|\langle x\rangle^2P_N(1-\chi_{R_N})K u\|_{L^2_{t,x}}^2
 &\lesssim\varepsilon^2\|w_au\|_{L^2_tH^s_x}^2.
\end{align*}
The two frequency costs are $(a+1)\theta=\ell$ and
$1-\theta(\rho-a-2)<\ell$.
The free massive Dirac homogeneous and weighted retarded estimates are
\[
 \|e^{-i(t-r)B_0}f_N\|_{L^2_tL^\infty_x}
 \lesssim_\eta N^{1+\eta}\|f_N\|_2,
\]
\[
 \left\|\int_r^t e^{-i(t-\tau)B_0}F_N(\tau)\,d\tau
 \right\|_{L^2_tL^\infty_x}
 \lesssim_\eta N^{1+\eta}
       \|\langle x\rangle^2F_N\|_{L^2_{t,x}}.
\]
They follow by splitting the free Dirac flow into its two scalar
massive half-wave branches; the weighted retarded proof is Lemma~\ref{dirac:free}, using the annular case of Lemma~\ref{paper:whitney}.
Combining these estimates with the two square sums gives
\[
 \sum_NN^{2(s-1-\eta-\ell)}
      \|v_N\|_{L^2_tL^\infty_x}^2\lesssim\|f\|_{H^s}^2.
\]
The initial exterior terms satisfy the corresponding estimate by the
varying-family Littlewood--Paley lemma.  The exponent is positive, so
Cauchy--Schwarz sums the exterior pieces.  This proves
\eqref{dirac:cubic-global-bound}.  Uniformity in $r$ and Minkowski
prove the $L^1H^s$ Duhamel assertion.
\end{proof}

Taking $a>3/2$ sufficiently close to $3/2$ shows that a sufficient
range for the cubic resolution estimate is
\begin{equation}\label{dirac:cubic-range}
 \rho>\frac72,\qquad
 s>1+\frac{5}{2(\rho-1)}.
\end{equation}
The local wave parametrix improves the interior estimate, while the
weighted energy controls the sum over all time intervals.

\section{Nonlinear estimates and scattering}
\label{sec:nonlinear}
\subsection{Spatial cube bounds and the Yukawa multiplier}

We first record explicitly the estimates used in the frequency
summation. Index the cube multipliers at width $\mu$ by their scaled
lattice centers $k\in\mathbb Z^3$, and denote them temporarily by
$C_k^\mu$. Their symbol supports have uniformly bounded overlap, and
\begin{equation}\label{tdnl:cubeframe}
 \sum_k\|C_k^\mu f\|_2^2\lesssim\|f\|_2^2,
 \qquad \sup_{\mu,k}\|\mathcal F^{-1}C_k^\mu\|_1<\infty.
\end{equation}
In the second expression the inverse Fourier transform is applied to
the multiplier symbol.

\begin{lemma}\label{tdnl:cubes}
Let $R_\mu$ have symbol supported in $\{|\xi|\leq C\mu\}$ and
a uniformly bounded $L^1$ convolution kernel. Then
\begin{align}
 \|R_\mu(f_L\overline g)\|_2
 &\lesssim\Gamma_{\mu,L}(f)\|g\|_2,
 \label{tdnl:cubeL2}\\
 \|R_\mu(f_L\overline{g_M})\|_\infty
 &\lesssim\Gamma_{\mu,L}(f)\Gamma_{\mu,M}(g).
 \label{tdnl:cubeLinf}
\end{align}
The same assertions hold for products without a conjugation.
\end{lemma}

\begin{proof}
Expand both factors in the cube partition. A term in a conjugate
product survives $R_\mu$ only when $|k-\ell|\leq C_0$; each index
therefore has a bounded number of partners. The triangle inequality
and Cauchy--Schwarz give
\[
 \begin{split}
 \|R_\mu(f_L\overline g)\|_2
 &\lesssim\sum_{|k-\ell|\leq C_0}
       \|C_k^\mu f_L\|_\infty\|C_\ell^\mu g\|_2\\
 &\lesssim
 \left(\sum_k\|C_k^\mu f_L\|_\infty^2\right)^{1/2}
 \left(\sum_\ell\|C_\ell^\mu g\|_2^2\right)^{1/2}.
 \end{split}
\]
Use \eqref{tdnl:cubeframe}. Putting both cube factors in $L^\infty$
proves \eqref{tdnl:cubeLinf}. Without a conjugation the surviving
indices satisfy $|k+\ell|\leq C_0$, and the proof is unchanged.
Every Cauchy--Schwarz inequality here is at a fixed time.
\end{proof}

Choose real even enlarged multipliers $\widetilde P_\mu$ equal to
one on the symbol support of $P_\mu$, with support in a fixed
enlargement of that block. Set
\[
 T_\mu F=V_b*P_\mu F,
 \qquad V_b(x)=\frac{e^{-b|x|}}{|x|},\quad b>0.
\]
With the Fourier convention $\widehat V_b(\xi)=4\pi/(b^2+|\xi|^2)$,
\begin{equation}\label{tdnl:yukawa}
 \|T_\mu\|_{L^2\to L^2}
 +\|T_\mu\|_{L^\infty\to L^\infty}
 \lesssim_b\mu^{-2}.
\end{equation}
For $\mu\geq2$, scale $\xi=\mu\eta$ and differentiate the symbol
\[
 4\pi\mu^{-2}\frac{p(\eta)}{(b/\mu)^2+|\eta|^2}.
\]
The cutoff $p$ is supported away from zero; its finitely many
derivatives needed for the $L^1$ kernel bound are uniform in
$b/\mu\geq0$. The $L^2$ bound follows directly from the multiplier
supremum. At $\mu=1$, use $\|V_b\|_1=4\pi b^{-2}$.
We shall use the identities
$T_\mu=T_\mu\widetilde P_\mu=
\widetilde P_\mu T_\mu$; no smooth cutoff is assumed idempotent.

\subsection{Summation in the physical space-time norm}

\begin{proposition}\label{tdnl:trilinear}
For $\mathcal T(u_1,u_2,u_3)=(V_b*(u_1\overline{u_2}))u_3$,
\begin{equation}\label{tdnl:trilinearbound}
 \|\mathcal T(u_1,u_2,u_3)\|_{L^1(I;H^s)}
 \lesssim_{s,b,\delta}
 \sum_{\pi\in S_3}
 Q_s(u_{\pi(1)};I)Q_s(u_{\pi(2)};I)
 \|u_{\pi(3)}\|_{L^\infty(I;H^s)}.
\end{equation}
In particular, if $\mathcal N(u)=\mathcal T(u,u,u)$, then
\begin{equation}\label{tdnl:difference}
 \|\mathcal N(u)-\mathcal N(v)\|_{L^1(I;H^s)}
 \lesssim_{s,b,\delta}
 (\|u\|_{W^s(I)}+\|v\|_{W^s(I)})^2
 \|u-v\|_{W^s(I)}.
\end{equation}
\end{proposition}

\begin{proof}
Begin with inputs having finitely many dyadic frequency blocks.
Such $L^2$ inputs are bounded continuous functions, and the cube
expansion at each fixed width has only finitely many nonzero terms.
Thus the following spatial pairings are justified directly.

Fix a time. Arrange the three input frequencies as $L\leq M\leq H$,
and let $\pi$ specify their placement among the three inputs. For an
$L^2$ test function $v_N$ in the output block, we claim
\begin{equation}\label{tdnl:quartic}
 \begin{split}
 &\left|\int T_\mu(u_{1,N_1}\overline{u_{2,N_2}})
                  u_{3,N_3}\overline{v_N}\,dx\right|\\
 &\quad\lesssim_b\mu^{-2}
       \Gamma_{\mu,L}(u_{\pi(1)})
       \Gamma_{\mu,M}(u_{\pi(2)})
       \|u_{\pi(3),H}\|_2\|v_N\|_2.
 \end{split}
\end{equation}
If $H=N_3$, apply \eqref{tdnl:cubeLinf} to the density, then
\eqref{tdnl:yukawa}, and put $u_{3,H}$ and $v_N$ in $L^2$.
If $H=N_1$, insert $\widetilde P_\mu$ into the density and into
the other product in the pairing. Estimate the former by
$\Gamma_{\mu,N_2}(u_2)\|u_{1,H}\|_2$, using
\eqref{tdnl:cubeL2}, and the latter by
$\Gamma_{\mu,N_3}(u_3)\|v_N\|_2$.
The $L^2$ bound in \eqref{tdnl:yukawa} completes this case.
If $H=N_2$, the same two $L^2$ estimates use $u_{1,N_1}$ and
$u_{3,N_3}$ as the cube-$L^\infty$ factors. The real even enlarged
cutoffs commute with conjugation, so these insertions also justify
the conjugations in each pairing. These three cases prove
\eqref{tdnl:quartic}.

Sum the interaction frequency first:
\begin{equation}\label{tdnl:musum}
 \sum_{\mu\in\mathcal D}\mu^{-2}
       \Gamma_{\mu,L}(u_i)\Gamma_{\mu,M}(u_j)
 \leq B_L(u_i)B_M(u_j).
\end{equation}
This is Cauchy--Schwarz with the same weight $\mu^{-2}$ in both
factors. Define the nonnegative sequences
\[
 a_N^j(t)=N^{s-\alpha}B_N(u_j)(t),\qquad
 e_N^j(t)=N^s\|u_{j,N}(t)\|_2.
\]
The cubic output has frequency at most a fixed multiple of $H$.
For a fixed $\kappa\geq1$, output duality and
\eqref{tdnl:quartic}--\eqref{tdnl:musum} therefore yield
\begin{equation}\label{tdnl:outputsum}
 \begin{split}
 N^s\|P_N\mathcal T(u_1,u_2,u_3)(t)\|_2
 &\lesssim_b\sum_{\pi\in S_3}
  \sum_{H\geq N/\kappa}(N/H)^s e_H^{\pi(3)}(t)\\
 &\quad\times\sum_{M\leq H}M^{-\delta}a_M^{\pi(2)}(t)
             \sum_{L\leq M}L^{-\delta}a_L^{\pi(1)}(t).
 \end{split}
\end{equation}
Separate applications of Cauchy--Schwarz show that the two inner
sums are at most
\[
 \left(\sum_{N\in\mathcal D}N^{-2\delta}\right)
       \|a^{\pi(1)}(t)\|_{\ell^2}\|a^{\pi(2)}(t)\|_{\ell^2}
 =\frac{\|a^{\pi(1)}(t)\|_{\ell^2}
             \|a^{\pi(2)}(t)\|_{\ell^2}}
          {1-2^{-2\delta}}.
\]
On logarithmic dyadic indices the remaining sum in $H$ is convolution
with $2^{-js}$, together with finitely many negative indices arising
from $\kappa$. Since $s>0$, this kernel belongs to $\ell^1$.
Discrete Young's inequality and Littlewood--Paley equivalence give
\begin{equation}\label{tdnl:pointwiseHs}
 \|\mathcal T(u_1,u_2,u_3)(t)\|_{H^s}
 \lesssim_{s,b,\delta}\sum_{\pi\in S_3}
   \|a^{\pi(1)}(t)\|_{\ell^2}
   \|a^{\pi(2)}(t)\|_{\ell^2}
   \|u_{\pi(3)}(t)\|_{H^s}.
\end{equation}
Integrate in time and apply Cauchy--Schwarz to the two $a$ factors.
This proves \eqref{tdnl:trilinearbound} for finite-band inputs.
In particular, the physical $H^s$ norm has been summed before time
integration; an $L_t^1$ norm has not been interchanged with an
output-frequency $\ell^2$ norm.

We now justify the passage to general inputs without a density
assumption on the global resolution norm. Let
$T_K=\chi_0(2^{-K}D)$ be smooth low-pass cutoffs. Their kernels have
uniformly bounded $L^1$ norms, and they commute with all the Fourier
projections above. Consequently
\begin{equation}\label{tdnl:truncationbounded}
 Q_s(T_Ku;I)\lesssim Q_s(u;I),\qquad
 \|T_Ku\|_{L^\infty H^s}\lesssim\|u\|_{L^\infty H^s}.
\end{equation}
For $s>\frac12$, choose a finite $p>3$ such that
$H^s(\mathbb R^3)\hookrightarrow L^p$, and put $r=p/(p-2)<3$.
The screened kernel belongs to $L^r$, so Young's inequality gives
\begin{equation}\label{tdnl:L2continuity}
 \|\mathcal T(f,g,h)\|_2
 \leq\|V_b\|_r\|f\|_p\|g\|_p\|h\|_2.
\end{equation}
Since $T_Ku(t)\to u(t)$ in $H^s$ at each time,
\eqref{tdnl:L2continuity} identifies the limit of the truncated
cubics in $L^2$. The $H^s$ norm is lower semicontinuous under strong
$L^2$ convergence. Fatou's lemma in time, together with
\eqref{tdnl:truncationbounded}, proves
\eqref{tdnl:trilinearbound} for the full inputs.
The resulting nonlinearity is strongly measurable in $H^s$:
the truncated expressions identify its weakly measurable $H^s$
limit almost everywhere, and $H^s$ is separable.
Finally, expand the difference of the two cubics into three polarized
terms to obtain \eqref{tdnl:difference}.
\end{proof}

The preceding scalar polarized estimate applies componentwise to
\eqref{main:hartree}; summation over finitely many spinor indices
only changes its constant. In particular,
\[
 \|\mathcal N_H(\psi)\|_{L_t^1H^s}
 \lesssim Q_s(\psi)^2\|\psi\|_{L_t^\infty H^s}.
\]
The cube norm used here is for the original spinor in the fixed frame;
no density multiplier is inserted into the convolution.

\subsection{Local cubic estimate}
Put $Z^s=L_t^\infty H^s\cap L_t^2L_x^\infty$. The fractional
product estimate~\cite{GrafakosOh2014}, for $s>0$, gives
\[
 \|(u_1^\dagger C_1u_2)C_2u_3\|_{L_t^1H^s}
 \lesssim\sum_{j=1}^3\|u_j\|_{L_t^\infty H^s}
                  \prod_{k\ne j}\|u_k\|_{L_t^2L_x^\infty}.
\]
It follows by expanding the difference into three polarized terms that
\[
 \|\mathcal N_C(u)-\mathcal N_C(v)\|_{L_t^1H^s}
 \lesssim(\|u\|_{Z^s}+\|v\|_{Z^s})^2\|u-v\|_{Z^s}.
\]
Neither this estimate nor the Hartree estimate uses a null structure.

\subsection{Exact-flow Duhamel, contraction, and asymptotic profiles}
Let $X^s=Z^s$ in the cubic case and $X^s=W^s$ in the Hartree case.
The homogeneous estimates proved above are uniform in the initial
time and hold in both orientations. For $G\in L_t^1H^s$, Minkowski
and restriction to time half-lines therefore give
\begin{equation}\label{dirac:duhamel}
 \left\|\int_0^tU_D(t,r)G(r)\,dr\right\|_{X^s}
 \lesssim\|G\|_{L_t^1H^s}.
\end{equation}
For instance, on $t\ge0$ write the integrand as
$\mathbf1_{t\ge r}U_D(t,r)G(r)$ with $r\ge0$; restrictions contract
each of the physical time norms. The past half-line is analogous.

\begin{lemma}[Uniqueness on intervals with small dispersive norm]
\label{paper:uniqueness}
Write $E(u;I)=\|u\|_{L^\infty(I;H^s)}$ and let
$D(u;I)=\|u\|_{L^2(I;L^\infty)}$ for the cubic equation, or
$D(u;I)=Q_s(u;I)$ for the Hartree equation. In either case the
polarized nonlinear bound gives, with $d=u-v$,
\begin{align*}
 \|\mathcal N(u)-\mathcal N(v)\|_{L^1(I;H^s)}
 \lesssim{}&[D(u;I)+D(v;I)]^2 E(d;I)\\
 &+[E(u;I)+E(v;I)][D(u;I)+D(v;I)]D(d;I).
\end{align*}
Consequently solutions with identical data at one endpoint agree
on $I$ whenever the dispersive norms are sufficiently small in
terms of their finite energy bounds.
\end{lemma}
\begin{proof}
Expand the cubic difference into three polarized terms. Terms in
which $d$ is the energy factor give the first line; the others give
the second. Apply the exact-flow $L^1H^s$ Duhamel estimate to $d$
and absorb. On each compact time interval the energy bounds are
finite and the squared dispersive norms are absolutely continuous
under time restriction. A finite subdivision therefore proves
uniqueness on that interval, and then on the full lifespan.
\end{proof}

\begin{proof}[Proof of Theorems~\ref{main:theorem} and \ref{main:hartree-theorem}]
The equation is equivalent to
\[
 \Phi(\psi)(t)=U_D(t,0)\psi_0+
       i\lambda\int_0^tU_D(t,r)\mathcal N(\psi(r))\,dr.
\]
The linear and nonlinear bounds show
\[
 \|\Phi(\psi)\|_{X^s}\le C_0\|\psi_0\|_{H^s}+C_1\|\psi\|_{X^s}^3,
\]
\[
 \|\Phi(\psi)-\Phi(\varphi)\|_{X^s}
 \le C_1(\|\psi\|_{X^s}+\|\varphi\|_{X^s})^2
                         \|\psi-\varphi\|_{X^s}.
\]
On the ball of radius $2C_0\|\psi_0\|_{H^s}$ the map is a strict
contraction for sufficiently small data. The mild solution is
continuous in $H^s$ because its source belongs to $L_t^1H^s$.
The same estimates prove local Lipschitz dependence. Uniqueness in the full stated resolution class
follows from Lemma~\ref{paper:uniqueness}.

The inverse propagator is uniformly bounded on $H^s$, so
\begin{equation}\label{dirac:profiles}
 \psi_\pm=\psi_0+i\lambda\int_0^{\pm\infty}
             U_D(0,r)\mathcal N(\psi(r))\,dr
\end{equation}
converges absolutely in $H^s$. The cocycle identity gives
\[
 \|\psi(t)-U_D(t,0)\psi_+\|_{H^s}
 \lesssim\int_t^\infty\|\mathcal N(\psi(r))\|_{H^s}\,dr\to0;
\]
the past limit is identical with its orientation reversed.
Given small $\psi_+$, solve instead
\[
 \psi(t)=U_D(t,0)\psi_+-i\lambda\int_t^\infty
                 U_D(t,r)\mathcal N(\psi(r))\,dr
\]
by the same contraction on the future half-line and extend the
small value at $0$ by the initial-value construction. This proves
the future wave operator; reverse time for the past operator.

\end{proof}

\appendix
\section{Atomic spaces adapted to the full evolution}
\label{sec:atomic}
We use the atomic and variation spaces of~\cite{HHK2009,HHK2010}.
For a Hilbert space $\mathcal H$, a $U^2(\mathcal H)$ atom is a step
function $a=\sum_{k=1}^{K}\mathbf1_{[t_{k-1},t_k)}f_{k-1}$,
where $-\infty=t_0<\cdots<t_K=+\infty$, $f_0=0$, and
$\sum_{k=0}^{K-1}\|f_k\|_{\mathcal H}^2=1$. The atomic norm is the
infimum of $\sum_j|c_j|$ over representations $u=\sum_jc_ja_j$.
The $V^p$ norm is the supremum, over finite increasing partitions,
of the $\ell^p$ norm of the increments. Endpoint conventions are
those of the erratum~\cite{HHK2010}: the value used at $+\infty$
in the variation sum is assigned to be zero, independently of the
actual limit. We use right-continuous profiles vanishing at
$-\infty$, after extending by zero before the initial endpoint.
In particular, constant profiles on a half-line are allowed and
are not required to tend to zero at its infinite endpoint.

\begin{proposition}[Transfer to exact-flow atomic spaces]
\label{paper:atomic}
Let $X^s(I)$ denote either resolution norm in the theorem and suppose
$\sup_r\|U_D(\cdot,r)f\|_{X^s(\mathbb R)}\le C\|f\|_{H^s}$.
On an oriented half-line define
\[
 \|u\|_{U_D^2(H^s)}
 :=\|U_D(0,t)u(t)\|_{U^2(H^s)},\qquad
 \|u\|_{V_D^2(H^s)}
 :=\|U_D(0,t)u(t)\|_{V^2(H^s)},
\]
using zero extension before the initial endpoint in the oriented
time variable. Then $\|u\|_{X^s(I)}\lesssim\|u\|_{U_D^2(H^s)}$.
If a solution has $\mathcal N(u)\in L^1(I;H^s)$, its pulled-back
profile belongs to both atomic and variation spaces, with bound
\[
 \|u\|_{U_D^2(H^s)}+\|u\|_{V_D^2(H^s)}
 \lesssim\|u(t_0)\|_{H^s}+|\lambda|\|\mathcal N(u)\|_{L^1(I;H^s)}.
\]
\end{proposition}
\begin{proof}
For a $U^2$ atom $a(t)=\sum_k\mathbf1_{I_k}(t)f_k$,
$\sum_k\|f_k\|_{H^s}^2=1$, the intervals $I_k$ are disjoint.
The dispersive part of $\|U_D(t,0)a(t)\|_{X^s}^2$ is therefore
the sum of its interval contributions, bounded by
$C\sum_k\|f_k\|_{H^s}^2$. The energy part is bounded by
$C\sup_k\|f_k\|_{H^s}$. Atomic decomposition proves the first
assertion. The profile identity and the uniform inverse propagator
bound give an absolutely continuous $H^s$-valued profile with total
variation bounded by the displayed right hand side, including its
jump at the initial endpoint. The standard embedding of right
continuous $V^1$ into $U^2$, and $V^1$ into $V^2$, gives the last
assertion. No Fourier projection is commuted through $U_D$.
\end{proof}

This transfer uses the estimates already proved for the exact
propagator. It does not assert an endpoint dispersive estimate or
a multilinear $V^2$ estimate beyond those established above.
The scattering profiles are defined relative to $U_D(t,r)$; a
comparison with a fixed Minkowski group would require a separate
linear scattering theorem.


\begin{thebibliography}{99}
\bibitem{MNO2003}
S. Machihara, K. Nakanishi, and T. Ozawa,
\emph{Small global solutions and the nonrelativistic limit for the nonlinear
Dirac equation}, Rev. Mat. Iberoam. \textbf{19} (2003), no.~1, 179--194.
\href{https://doi.org/10.4171/RMI/342}{doi:10.4171/RMI/342}.

\bibitem{BH2015}
I. Bejenaru and S. Herr,
\emph{The cubic Dirac equation: small initial data in
$H^1(\mathbb R^3)$}, Comm. Math. Phys. \textbf{335} (2015), no.~1, 43--82.
\href{https://doi.org/10.1007/s00220-014-2164-0}{doi:10.1007/s00220-014-2164-0}.

\bibitem{HT}
S. Herr and A. Tesfahun,
\emph{Small data scattering for semi-relativistic equations with Hartree type
nonlinearity}, J. Differential Equations \textbf{259} (2015), no.~10,
5510--5532.
\href{https://doi.org/10.1016/j.jde.2015.06.037}{doi:10.1016/j.jde.2015.06.037}.

\bibitem{Yang2019}
C. Yang,
\emph{Scattering results for Dirac Hartree-type equations with small initial
data}, Commun. Pure Appl. Anal. \textbf{18} (2019), no.~4, 1711--1734.
\href{https://doi.org/10.3934/cpaa.2019081}{doi:10.3934/cpaa.2019081}.

\bibitem{Tesfahun2020}
A. Tesfahun,
\emph{Small data scattering for cubic Dirac equation with Hartree type
nonlinearity in $\mathbb R^{1+3}$}, SIAM J. Math. Anal. \textbf{52} (2020),
no.~3, 2969--3003.
\href{https://doi.org/10.1137/17M1155788}{doi:10.1137/17M1155788}.

\bibitem{CHL2023}
Y. Cho, S. Hong, and K. Lee,
\emph{Scattering and nonscattering of the Hartree-type nonlinear Dirac system
at critical regularity}, SIAM J. Math. Anal. \textbf{55} (2023), no.~4,
3395--3419.
\href{https://doi.org/10.1137/21M1437135}{doi:10.1137/21M1437135}.

\bibitem{DF2008}
P. D'Ancona and L. Fanelli,
\emph{Strichartz and smoothing estimates for dispersive equations with
magnetic potentials}, Comm. Partial Differential Equations \textbf{33}
(2008), 1082--1112.
\href{https://doi.org/10.1080/03605300701743749}{doi:10.1080/03605300701743749}.

\bibitem{CdS2019}
F. Cacciafesta and A.-S. de Suzzoni,
\emph{Weak dispersion for the Dirac equation on asymptotically flat and
warped product spaces}, Discrete Contin. Dyn. Syst. \textbf{39} (2019),
no.~8, 4359--4398.
\href{https://doi.org/10.3934/dcds.2019177}{doi:10.3934/dcds.2019177}.

\bibitem{CdSM2024}
F. Cacciafesta, A.-S. de Suzzoni, and L. Meng,
\emph{Strichartz estimates for the Dirac equation on asymptotically flat
manifolds}, Ann. Sc. Norm. Super. Pisa Cl. Sci. (5) \textbf{25} (2024),
no.~3, 1799--1812.
\href{https://doi.org/10.2422/2036-2145.202203_026}{doi:10.2422/2036-2145.202203\_026}.

\bibitem{Smith1998}
H. F. Smith,
\emph{A parametrix construction for wave equations with $C^{1,1}$
coefficients}, Ann. Inst. Fourier (Grenoble) \textbf{48} (1998), no.~3,
797--835.
\href{https://doi.org/10.5802/aif.1640}{doi:10.5802/aif.1640}.

\bibitem{MT2012}
J. Metcalfe and D. Tataru,
\emph{Global parametrices and dispersive estimates for variable coefficient
wave equations}, Math. Ann. \textbf{353} (2012), no.~4, 1183--1237.
\href{https://doi.org/10.1007/s00208-011-0714-8}{doi:10.1007/s00208-011-0714-8}.

\bibitem{HH2026}
S. Herr and S. Hong,
\emph{Strichartz estimates for the half Klein--Gordon equation on
asymptotically flat backgrounds and applications to cubic Dirac equations},
arXiv:2502.13670v3 [math.AP], 2026 (first submitted 2025).
\url{https://arxiv.org/abs/2502.13670v3}.

\bibitem{HHK2009}
M. Hadac, S. Herr, and H. Koch,
\emph{Well-posedness and scattering for the KP-II equation in a critical
space}, Ann. Inst. H. Poincar\'e Anal. Non Lin\'eaire \textbf{26} (2009),
no.~3, 917--941.
\href{https://doi.org/10.1016/j.anihpc.2008.04.002}{doi:10.1016/j.anihpc.2008.04.002}.

\bibitem{HHK2010}
M. Hadac, S. Herr, and H. Koch,
\emph{Erratum to ``Well-posedness and scattering for the KP-II equation in
a critical space''}, Ann. Inst. H. Poincar\'e Anal. Non Lin\'eaire
\textbf{27} (2010), no.~3, 971--972.
\href{https://doi.org/10.1016/j.anihpc.2010.01.006}{doi:10.1016/j.anihpc.2010.01.006}.


\bibitem{Shubin2001}
M. A. Shubin, \emph{Pseudodifferential Operators and Spectral Theory},
second edition, Springer, Berlin, 2001.
\href{https://doi.org/10.1007/978-3-642-56579-3}{doi:10.1007/978-3-642-56579-3}.

\bibitem{GrafakosOh2014}
L. Grafakos and S. Oh, \emph{The Kato--Ponce inequality},
Comm. Partial Differential Equations \textbf{39} (2014), no.~6,
1128--1157.
\href{https://doi.org/10.1080/03605302.2013.822885}{doi:10.1080/03605302.2013.822885}.

\bibitem{Alcubierre}
M.~Alcubierre, \emph{The Dirac equation in general relativity and the
3+1 formalism}, General Relativity and Gravitation \textbf{57} (2025),
article 122.
\href{https://doi.org/10.1007/s10714-025-03460-8}{doi:10.1007/s10714-025-03460-8}.
\end{thebibliography}
\end{document}